\documentclass[11pt, reqno]{amsart}

\numberwithin{equation}{section}

\usepackage{amscd,amsmath}
\usepackage{mathrsfs}
\usepackage{amsfonts}
\usepackage{amssymb}
\usepackage{enumerate}
\usepackage{setspace}
\usepackage{color}

\usepackage[margin=1in]{geometry}

\usepackage[colorlinks=true, citecolor=blue]{hyperref}

\newtheorem{theorem}{Theorem}[section]

\newtheorem{remark}{Remark}[section]

\newcommand{\eqn}{\begin{eqnarray}}
\newcommand{\een}{\end{eqnarray}}

\newcommand{\la}{\langle}
\newcommand{\ra}{\rangle}

\DeclareMathOperator{\dv}{div}
\DeclareMathOperator{\curl}{curl}

\begin{document}

\title[$2\frac{1}{2}$ dimensional Hall equation]{On the existence and temporal asymptotics  of solutions for the two and half dimensional Hall MHD}

\author{Hantaek Bae}
\address{Department of Mathematical Sciences, Ulsan National Institute of Science and Technology, Korea}
\email{hantaek@unist.ac.kr}

\author{Kyungkeun Kang}
\address{Department of Mathematics, Yonsei University, Korea}
\email{kkang@yonsei.ac.kr}

\date{\today}

\begin{abstract}
In this paper, we deal with the $2\frac{1}{2}$ dimensional Hall MHD by taking the velocity field $u$ and the magnetic field $B$ of the form $u(t,x,y)=\left(\nabla^{\perp}\phi(t,x,y), W(t,x,y)\right)$ and $B(t,x,y)=\left(\nabla^{\perp}\psi(t,x,y), Z(t,x,y)\right)$. 

We begin with the Hall equations (without the effect of the fluid part). We first show the long time behavior of weak solutions and weak-strong uniqueness. We then proceed to prove the existence of unique strong solutions locally in time and to derive a blow-up criterion. We also demonstrate that the strong solution exists globally in time and decay algebraically  if some smallness conditions are imposed. We further improve the decay rates of $\psi$ using the structure of the equation of $\psi$. As a consequence of the decay rates of $(\psi,Z)$, we find the asymptotic profiles of $(\psi,Z)$. We finally show that a small perturbation of initial data near zero can be extended to a small perturbations near harmonic functions. 

In the presence of the fluid filed, the results, by comparison, fall short of the previous ones in the absence of the fluid part. We prove two results: the existence of unique strong solutions locally in time and a blow-up criterion, and  the existence of unique strong solutions globally in time with some smallness condition on initial data.
\end{abstract}

\maketitle

\vspace{-2ex}

\section{Introduction}
The Magnetohydrodynamics equations (MHD in short) provide a macroscopic description of a plasma and provide a relevant description for fusion plasmas, the solar interior and its atmosphere, the Earth's magnetosphere and inner core, etc. The governing equations for the incompressible and resistive MHD are
\begin{subequations} \label{MHD}
\begin{align}
\text{Momentum Equation:} \quad &  u_{t}  +u\cdot \nabla u - J\times B+\nabla p-\mu\Delta u=0, \label{MHD a}\\
\text{Incompressibility:} \quad & \dv u=0, \label{MHD b}\\
\text{Amp$\grave{\text{e}}$re's Law:} \quad & \curl B=\mu_{0}J, \label{MHD c}\\
\text{Faraday's Law:} \quad & \curl E=-B_{t},\label{MHD d}\\
\text{Ohm's Law for resistive MHD:} \quad &E+u\times B=\nu J, \label{MHD e}\\
\text{Incompressibility:} \quad& \dv B=0, \label{MHD f}
\end{align}
\end{subequations}
where $u$ is the  velocity field, $p$ is the pressure, and $B$ is the magnetic field. $\mu$ and $\nu$ are the viscosity and the resistivity constants, respectively. The right-hand side of (\ref{MHD e}) is called the collision term and $J\times B$ in (\ref{MHD a}) is called the Lorentz force. However, (\ref{MHD}) is deficient in many respect: for example, (\ref{MHD}) does not explain magnetic reconnection on the Sun which is very important role in acceleration plasma by converting magnetic energy into bulk kinetic energy. For this reason, the generalized Ohm's Law is required and we here take the following
\eqn \label{Generalized Ohm}
E+u\times B=\nu J+\frac{1}{en}\left(J\times B-\nabla p_{e}\right),
\een 
where  $e$ is the elementary charge, $n$ is the number density, and $p_{e}$ is the electron pressure. The  second term on the right-hand side of (\ref{Generalized Ohm})  is called the Hall term. In terms of $(u,\overline{p},B)$, we have the Hall MHD with $\mu_{0}=en=1$ for simplicity:
\begin{subequations}\label{Hall MHD}
\begin{align}
&u_{t}  +u\cdot \nabla u - B\cdot \nabla B+\nabla \overline{p}-\mu\Delta u=0,  \label{Hall MHD a}\\
&B_{t}  + u\cdot \nabla B -B\cdot \nabla u +\curl \left((\curl B)\times B\right)- \nu\Delta B=0,  \\
&\dv u=0, \quad \dv B=0,
\end{align}
\end{subequations}
where we use the following in (\ref{Hall MHD a}):
\[
J\times B=B\cdot \nabla B-\frac{1}{2}\nabla |B|^{2}, \quad \overline{p}=p+\frac{1}{2}|B|^{2}.
\]

The Hall-MHD is important in describing many physical phenomena \cite{Balbus, Forbes, Homann, Lighthill, Mininni, Shalybkov, Wardle}. The Hall-MHD recently  has been studied intensively.  The Hall-MHD can be derived from either two fluids model or kinetic models in a mathematically rigorous way  \cite{Acheritogaray}.  Global weak solution, local classical solution, global solution for small data, and decay rates are established in \cite{Chae-Degond-Liu, Chae-Lee, Chae Schonbek}. There have been many follow-up results of these papers; see \cite{Chae Wan Wu, Chae Weng, Dai, Fan, Han, Wan 1, Wan 2, Wan 3, Wan 4, Weng 1, Weng 2} and references therein.

\subsection{$2\frac{1}{2}$ dimensional Hall MHD}
The Hall term, $\curl \left(\left(\curl B\right)\times B\right)$, is dominant  when analyzing (\ref{Hall MHD}). So, even if we deal with (\ref{Hall MHD}) in the $2 \frac{1}{2}$ dimensional case, the global regularity problem for  (\ref{Hall MHD})  is still open. As a result, compared to MHD  and the incompressible Navier-Stokes equations, there are only a few results dealing with the $2 \frac{1}{2}$ dimensional (\ref{Hall MHD}); partial regularity theory \cite{Chae Wolf 1}, global regularity with partial dissipations in (\ref{Hall MHD a}) \cite{Du}, irreducibility property \cite{Yamazaki}.

In this paper,  we take its $2 \frac{1}{2}$ dimensional form of (\ref{Hall MHD}) through 
\begin{subequations}\label{2D reduction}
\begin{align}
u(t,x,y)&=\left(\nabla^{\perp}\phi(t,x,y), W(t,x,y)\right)=\left(-\phi_{y}(t,x,y), \phi_{x}(t,x,y), W(t,x,y)\right), \label{2D reduction a}\\
B(t,x,y)&=\left(\nabla^{\perp}\psi(t,x,y), Z(t,x,y)\right)=\left(-\psi_{y}(t,x,y), \psi_{x}(t,x,y), Z(t,x,y)\right). \label{2D reduction b}
\end{align}
\end{subequations}
Due to the presence of the pressure in (\ref{Hall MHD a}), we take the curl to (\ref{Hall MHD a}) and we rewrite (\ref{Hall MHD}) as 
\begin{subequations} \label{coupled two half}
\begin{align}
& \psi_{t}-\Delta \psi =[\psi,Z]-[\psi,\phi],\\
& Z_{t}-\Delta Z=[\Delta \psi,\psi]-[Z,\phi]+[W,\psi],\\
& W_{t}-\Delta W=-[W,\phi]-[\psi,Z],\\
& \Delta \phi_{t}-\Delta^{2} \phi=-[\Delta \phi,\phi]+[\Delta \psi,\psi],
\end{align}
\end{subequations}
where we set $\mu=\nu=1$ for simplicity and $[f,g]=\nabla f\cdot \nabla^{\perp}g=f_{x}g_{y}-f_{y}g_{x}$. (\ref{coupled two half})  is used to investigate the influence of the Hall-term on the island width of a tearing instability \cite{Homann} and to show a finite-time collapse to a current sheet \cite{Brizard, Janda 1, Janda 2, Litvinenko}. (\ref{coupled two half})  is also used in \cite{Chae Wolf 2} to study regularity of stationary weak solutions.

\subsection{\bf Hall equation}
Since the Hall term  is dominant when we deal with (\ref{coupled two half}), we will mainly concentrate on the Hall equations: (\ref{coupled two half}) without the effect of the fluid part. Then, (\ref{coupled two half}) is reduced to the following equations
\begin{subequations}\label{Hall equation in 2D}
\begin{align}
&\psi_{t}-\Delta \psi=[\psi,Z], \label{Hall equation in 2D a}\\
&Z_{t} -\Delta Z=[\Delta \psi,\psi].   \label{Hall equation in 2D b}
\end{align}
\end{subequations}
We will explain several results of (\ref{Hall equation in 2D}) from Section \ref{sec:1.2} to Section \ref{sec:1.5}, but before we do, we briefly describe them.
\begin{enumerate}[]
\item (1) The existence of weak solutions and decay rates of (\ref{Hall equation in 2D}) can be proved by following \cite{Chae-Degond-Liu, Chae Schonbek}. In Section \ref{sec:1.2}, we restate the decay rate of weak solutions in \cite{Chae Schonbek} to (\ref{Hall equation in 2D}) (Theorem \ref{weak solution}) and establish weak-strong uniqueness (Theorem \ref{weak strong uniqueness}).
\item (2) In Section \ref{sec:1.3}, we deal with strong solutions of (\ref{Hall equation in 2D}). We first establish the existence of unique local-in-time solutions with large initial data and  a blow-up criterion (Theorem \ref{LWP}). Having established the local in time results, we then proceed to extend the solution globally-in-time and to derive decay rates by imposing some smallness condition to initial data (Theorem \ref{GWP}). We also improve the decay rates of $\psi$ by using the structure of the equation of $\psi$ (Theorem \ref{decay of psi}).
\item (3) It is reasonable to study the asymptotic stability of temporally decaying solutions in Theorem \ref{GWP} and Theorem \ref{decay of psi}. In Section \ref{sec:1.4}, we intend to find asymptotic profiles of such  solutions of (\ref{Hall equation in 2D}) from the observation that constant multiples of the two dimensional heat kernel $\Gamma$ are solutions of (\ref{Hall equation in 2D}) (Theorem \ref{Asymptotics}).
\item (4) The aim of Section \ref{sec:1.5} is to analyze (\ref{Hall equation in 2D}) around harmonic functions. We take $\psi=\rho+\overline{\psi}$ or $Z=\omega+\overline{Z}$, where $\overline{\psi}$ and $\overline{Z}$ are harmonic functions. We show that there exists unique global-in-time solutions if $\rho_{0}$ or $\omega_{0}$ are sufficiently small (Theorem \ref{perturbation theorem 1} and Theorem \ref{perturbation theorem 2}). We emphasize that the size of $\overline{\psi}$ and $\overline{Z}$ are arbitrary. 
\end{enumerate}

\subsection{\bf Weak solution of (\ref{Hall equation in 2D})}\label{sec:1.2}
The the existence and decay rate of a weak solution even for (\ref{Hall MHD}) are  already proved in \cite{Chae-Degond-Liu, Chae Schonbek} with $u_{0}\in L^{2}$ and $B_{0}\in L^{2}$. We here restate these results to (\ref{Hall equation in 2D}). We first note that  we can derive the following:
\[
\frac{1}{2}\frac{d}{dt} \left(\left\|\nabla \psi\right\|^{2}_{L^{2}}+\left\|Z\right\|^{2}_{L^{2}}\right) +\left\|\Delta \psi\right\|^{2}_{L^{2}}+ \left\|\nabla Z\right\|^{2}_{L^{2}}=0.
\]
This is enough to show the existence of a weak solution of (\ref{Hall equation in 2D}) with $\left(\nabla\psi_{0}, Z_{0}\right)\in L^{2}$. Moreover, we have temporal decay rates of weak solutions which is the two dimensional  version of \cite{Chae Schonbek}.

\begin{theorem}\label{weak solution}
Let $\left(\nabla\psi_{0}, Z_{0}\right)\in L^{2}$. Then, there is a weak solution of (\ref{Hall equation in 2D}) satisfying
\[
\left\|\nabla \psi(t)\right\|^{2}_{L^{2}}+\left\|Z(t)\right\|^{2}_{L^{2}}+2\int^{t}_{0}\left(\left\|\Delta \psi(s)\right\|^{2}_{L^{2}}+\left\|\nabla Z(s)\right\|^{2}_{L^{2}}\right)ds \leq \left\|\nabla \psi_{0}\right\|^{2}_{L^{2}}+\left\|Z_{0}\right\|^{2}_{L^{2}}
\]
for all $t>0$.  If $\left(\nabla\psi_{0}, Z_{0}\right)\in L^{2}\cap L^{1}$,  $\nabla\psi$ and $Z$ decay in time as
\eqn\label{Weak decay}
\left\|\nabla \psi(t)\right\|_{L^{2}}+\left\|Z(t)\right\|_{L^{2}}\leq \frac{C_{0}}{\sqrt{1+t}},
\een
where $C_{0}$ depends on $\left\|\nabla \psi_{0}\right\|_{L^{2}\cap L^{1}}$ and $\left\|Z_{0}\right\|_{L^{2}\cap L^{1}}$.
\end{theorem}

As in the case of the incompressible Navier-Stokes equations, the uniqueness of weak solution of (\ref{Hall equation in 2D}) is unknown. Weak-strong uniqueness is to find  a path space $\mathcal{P}$ of a strong solution $B\in \mathcal{P}$ such that all weak solutions which share the same initial condition $B_{0}$ equal $B$. In this paper, we do not aim to derive very general weak-strong uniqueness results as in the case of the incompressible Navier-Stokes equations \cite{Germain}, but focus on Serrin-type results.

\begin{theorem}\label{weak strong uniqueness}
Let $B_{1}=\left(\nabla^{\perp}\psi_{1}, Z_{1}\right)$ and $B_{2}=\left(\nabla^{\perp}\psi_{2}, Z_{2}\right)$ be weak solutions of (\ref{Hall equation in 2D})  with the same initial data $\left(\nabla\psi_{0}, Z_{0}\right)\in L^{2}$. Then, $B_{1}=B_{2}$ on $[0,T]$ if  $B_{2}$ satisfies the condition 
\eqn \label{Serrin condition}
\left(\Delta\psi_{2}, \nabla Z_{2}\right)\in L^{p}\left([0,T]: L^{q}\right), \quad \frac{1}{p}+\frac{1}{q}=\frac{1}{2}, \quad 2\leq q<\infty.
\een
\end{theorem}

\begin{remark} \upshape
In \cite{Chae-Degond-Liu}, weak-strong uniqueness is established with $B_{2}\in L^{2}\left([0,T]: W^{1, \infty}(\mathbb{R}^{3})\right)$. By contrast, we derive weak-strong uniqueness with $B_{2}\in L^{p}\left([0,T]: L^{q}(\mathbb{R}^{2})\right)$.
\end{remark}

\subsection{\bf Strong solutions of (\ref{Hall equation in 2D})} \label{sec:1.3}
In this paper, we establish the local  in time existence of unique strong solutions of (\ref{Hall equation in 2D}) with initial data $\left(\nabla\psi_{0}, Z_{0}\right)\in H^{2}$. Let
\begin{equation}  \label{energy norm}
\begin{split}
&M(t)=\left\|\nabla \psi(t)\right\|^{2}_{H^{2}}+\left\|Z(t)\right\|^{2}_{H^{2}}, \quad M(0)=\left\|\nabla \psi_{0}\right\|^{2}_{H^{2}}+\left\|Z_{0}\right\|^{2},\\
&N(t)=\left\|\nabla^{2} \psi(t)\right\|^{2}_{H^{2}}+\left\|\nabla Z(t)\right\|^{2}_{H^{2}}, \quad \mathcal{E}(t)=M(t)+\int^{t}_{0}N(s)ds.
\end{split}
\end{equation} 
When the energy method is applied, we observe that  the terms with the highest derivative, that are unlikely to be handled by Laplacian's regularity, disappear due to the the properties of the commutator in (\ref{commutator}). For example, see (\ref{H2 bound 1}). So, we can derive the following inequalities: 
\begin{equation*}
\begin{split}
&\frac{d}{dt}(1+M)+N \leq C\left(1+M\right)^{3} 
\end{split}
\end{equation*}
and this gives the first part of the following result. To derive a blow-up criterion, we re-estimates some terms in Section \ref{sec:3.1.1} from $L^{2}$ to $L^{r}$, $r\ne 2$.  From now on, constants that depend on $M(0)$ are not specifically specified each time when we state our results, and we will use $\mathcal{E}_{0}$ in common.

\begin{theorem} \label{LWP}
Let $\left(\nabla\psi_{0}, Z_{0}\right)\in H^{2}$. There exists $T^{\ast}=T(\mathcal{E}_{0})$ such that there exists a unique solution  of (\ref{Hall equation in 2D}) with $\mathcal{E}(T^{\ast})<\infty$. Moreover, the maximal existence time $T^{\ast}<\infty$ if and only if 
\eqn \label{blowup}
\lim_{T\nearrow T^{\ast}}\int^{T}_{0}\left\|\nabla Z(t)\right\|^{q}_{L^{p}}dt=\infty, \quad \frac{1}{p}+\frac{1}{q}=\frac{1}{2}, \quad 2\leq q<\infty.
\een
\end{theorem}

\begin{remark} \upshape
Compared to \cite{Chae-Degond-Liu}, the regularity of initial data  is the borderline case: $B_{0}=\left(\nabla^{\perp}\psi_{0}, Z_{0}\right)\in H^{2}$ with $2=\frac{d}{2}+1=\frac{2}{2}+1$. Moreover, the  blow-up criterion in (\ref{blowup})  is stated only in terms of the third component of $B$.  A similar blow-up criterion can be derived in terms of  $\Delta\psi$: 
\eqn  \label{blowup 2222}
T^{\ast}<\infty \iff \lim_{T\nearrow T^{\ast}}\int^{T}_{0}\left\|\Delta \psi\right\|^{q}_{L^{p}}dt=\infty, \quad \frac{1}{p}+\frac{1}{q}=\frac{1}{2}, \quad 2\leq q<\infty.
\een
The condition (\ref{Serrin condition}) and the blow-up criteria (\ref{blowup}) and (\ref{blowup 2222})  are related to the scaling invariant property of  (\ref{Hall equation in 2D}): if $(\psi(t,x,y),Z(t,x,y))$ is a solution of (\ref{Hall equation in 2D}) on $[0,T]$, so is 
\eqn \label{scaling invariance}
(\psi_{\lambda}(t,x,y)=\lambda^{-1}\psi\left(\lambda^{2}t,\lambda x, \lambda y\right), \quad Z_{\lambda}(t,x,y)=Z\left(\lambda^{2}t,\lambda x, \lambda y\right) \quad \text{on $[0,\lambda^{2}T]$}.
\een  
\end{remark}
\vspace{1ex}

Since (\ref{Hall equation in 2D}) is dissipative, we typically expect the global well-posedness and temporal decay rates of solutions under some smallness conditions. In Section \ref{sec:3}, we will derive the followings
\begin{equation} \label{two inequalities Hall}
\begin{split}
&\frac{d}{dt}\left(\left\|\Delta \psi\right\|^{2}_{L^{2}}+\left\|\nabla Z\right\|^{2}_{L^{2}}\right)+\left\|\nabla \Delta \psi\right\|^{2}_{L^{2}}+\left\|\Delta Z\right\|^{2}_{L^{2}}\leq CS(t)\left(\left\|\nabla \Delta \psi\right\|^{2}_{L^{2}}+\left\|\Delta Z\right\|^{2}_{L^{2}}\right),\\
&\frac{d}{dt}\left(\left\|\nabla \Delta \psi\right\|^{2}_{L^{2}}+\left\|\Delta Z\right\|^{2}_{L^{2}}\right)+\left\|\Delta^{2} \psi\right\|^{2}_{L^{2}}+\left\|\nabla \Delta Z\right\|^{2}_{L^{2}}\leq CS(t)\left(\left\|\Delta^{2} \psi\right\|^{2}_{L^{2}}+\left\|\nabla \Delta Z\right\|^{2}_{L^{2}}\right),
\end{split}
\end{equation}
where $S(t)=\left\|\Delta \psi\right\|^{2}_{L^{2}}+\left\|\nabla Z\right\|^{2}_{L^{2}}$. By imposing the smallness condition of the form
\eqn \label{smallness condition 1}
\epsilon_{1}=\left\|\Delta \psi_{0}\right\|^{2}_{L^{2}}+ \left\|\nabla Z_{0}\right\|^{2}_{L^{2}}, \quad C\epsilon_{1}<1
\een
we can obtain global-in-time solutions and can find decay rates of the solution in Theorem \ref{LWP}.

\begin{theorem} \label{GWP}
Let $(\nabla \psi_{0}, Z_{0})\in H^{2}$ which satisfies (\ref{smallness condition 1}). Then, we can take $T^{\ast}=\infty$ in Theorem \ref{LWP}. If $(\nabla \psi_{0}, Z_{0})\in L^{1}$ in addition, $(\Delta \psi, \nabla Z)$ decays in time as follows
\eqn \label{GWP Decay 2}
\left\|\Delta \psi(t)\right\|_{L^{2}}+\left\|\nabla Z(t)\right\|_{L^{2}} \leq \frac{\mathcal{E}_{0}}{1+t},\quad \left\|\nabla \Delta \psi(t)\right\|_{L^{2}}+\left\|\Delta Z(t)\right\|_{L^{2}}\leq \frac{\mathcal{E}_{0}}{(1+t)^{3/2}}.
\een
\end{theorem}

\begin{remark} \upshape
\begin{enumerate}[] 
\item (1) The more detailed dependence of $\mathcal{E}_{0}$  on the initial condition is given in Section \ref{sec:3.3.2}. From the linear part of (\ref{Hall equation in 2D}), we expect the decay rates of the form  
\[
\left\|\Delta \psi(t)\right\|_{L^{2}}+\left\|\nabla Z(t)\right\|_{L^{2}} \leq \frac{\mathcal{E}_{0}}{1+\sqrt{t}},\quad \left\|\nabla \Delta \psi(t)\right\|_{L^{2}}+\left\|\Delta Z(t)\right\|_{L^{2}}\leq \frac{\mathcal{E}_{0}}{1+t}
\]
but, these are improved to (\ref{GWP Decay 2})  by combining with (\ref{Weak decay}).  
\item (2) According to (\ref{scaling invariance}), $\dot{H}^{1}$ is a scaling-invariant space of $(\nabla\psi_{0}, Z_{0})$ and  we use the smallness condition (\ref{smallness condition 1}) in Theorem \ref{GWP}. In \cite{Chae-Lee}, the smallness condition is presented in the Besov space $\dot{B}^{3/2}_{2,1}$. We here replace $\dot{B}^{3/2}_{2,1}$ with $\dot{H}^{1}$. This is possible by two reasons: (i) we reduce the dimension from 3 to 2 and (ii) we are able to avoid the Littlewood-Paley theory by exploiting some cancellation properties of the commutator (\ref{commutator}).
\item (3) In the proof of Theorem \ref{LWP}, we obtain (\ref{H2 bound 2}):
\[
\frac{d}{dt}\left(\left\|\Delta \psi\right\|^{2}_{L^{2}}+\left\|\nabla Z\right\|^{2}_{L^{2}}\right)+\left\|\nabla \Delta \psi\right\|^{2}_{L^{2}}+\left\|\Delta Z\right\|^{2}_{L^{2}}\leq C\left\|\Delta \psi\right\|^{2}_{L^{2}} \left\|\nabla\Delta \psi\right\|^{2}_{L^{2}}.
\]
Combined with the uniqueness part in Section \ref{sec:3.1.2}, we can show the existence of a unique solution globally-in-time with $(\nabla \psi_{0}, Z_{0})\in \dot{H}^{1}$ under the smallness condition (\ref{smallness condition 1}). But, we do not state this case separately because we are more interested in using Theorem \ref{GWP} to prove Theorem \ref{decay of psi} and Theorem \ref{Asymptotics}.
\end{enumerate}
\end{remark}

The decay rates in Theorem \ref{GWP} are obtained by treating $\nabla \psi$ and $Z$ together. But, we observe  that we can improve the decay rates of $\psi$  by using the structure of (\ref{Hall equation in 2D a}) which is a dissipative transport equation, and  this is also the reason why the same method cannot be applied to $Z$. 

\begin{theorem} \label{decay of psi}
Let $(\nabla \psi_{0}, Z_{0})\in H^{2}$ satisfy (\ref{smallness condition 1}). If $\psi_{0}\in L^{1}\cap L^{2}$ in addition, $\psi$  decays in time as follows:
\eqn  \label{improved decay psi}
\left\|\nabla \psi(t)\right\|_{L^{2}}\leq \frac{\mathcal{E}_{0}}{1+t}, \quad \left\|\Delta \psi(t)\right\|_{L^{2}}\leq \frac{\mathcal{E}_{0}}{(1+t)^{3/2}}  \quad \text{for all $t>0$}.
\een
\end{theorem}

\subsection{Asymptotic behaviors} \label{sec:1.4}
Theorem \ref{GWP} and  Theorem \ref{decay of psi} provide upper bounds of decay rates. In particular, we have
\[
\left\|\nabla \psi(t)\right\|_{L^{2}}\leq \frac{\mathcal{E}_{0}}{1+t}, \quad \left\|\nabla Z(t)\right\|_{L^{2}}\leq \frac{\mathcal{E}_{0}}{1+t}.
\]
Although there is no embedding relationship between $\dot{H}^{1}$ and $L^{\infty}$, we can expect similar decay results in $L^{\infty}$ if we establish the asymptotic behavior of $(\psi,Z)$ as $t\rightarrow \infty$. One motivation for considering  asymptotic behavior is that the asymptotic behavior of the vorticity of the incompressible Navier-Stokes equations  in 2D is well-established: \cite{Carpio}, \cite{Gallay}, \cite[Page 44]{Giga}. And similar results can be obtained to more complicated models such as an aerotaxis model coupled to fluid equations \cite{Chae Kang Lee}.  Along this direction, we also want to find an asymptotic profile of the solutions of (\ref{Hall equation in 2D}). To do so, we assume $\left(\nabla\psi_{0}, Z_{0}\right)\in H^{2}$ as before and  we impose the following additional conditions
\begin{subequations}\label{L1 average}
\begin{align}
&\psi_{0}\in L^{1}, \quad Z_{0}\in L^{1}, \label{L1 average a}\\
& \la x\ra \psi_{0}\in L^{1},  \quad \la x\ra Z_{0}\in L^{1}, \label{L1 average b}\\
& \int_{\mathbb{R}^{2}}\psi_{0}(x)dx=\gamma, \quad \int_{\mathbb{R}^{2}}Z_{0}(x)dx=\eta, \label{L1 average c}
\end{align}
\end{subequations}
where $\la x \ra=\sqrt{1+|x|^{2}}$. Depending on which of one given by (\ref{L1 average}) we choose, we can describe asymptotic behaviors accordingly.

\begin{theorem}\label{Asymptotics}
Suppose $\left(\nabla\psi_{0}, Z_{0}\right)\in H^{2}$ and we assume (\ref{L1 average a}). Then, we obtain
\[
\psi(t,x)=\Gamma(t) \ast \psi_{0}+O(t^{-3/2}),\quad Z(t,x)=\Gamma(t) \ast Z_{0} +O(t^{-2})
\]
as $t\rightarrow \infty$, where $\Gamma$ is the two dimensional heat kernel. If we assume (\ref{L1 average b}) and (\ref{L1 average c}),
\[
\psi(t,x)=\gamma\Gamma(t,x) +O(t^{-3/2}),\quad Z(t,x)=\eta\Gamma(t,x) +O(t^{-3/2}).
\]
\end{theorem}

We note that the asymptotic behavior of $\nabla \psi$ is the same as $Z$ because $\nabla\left(\nabla^{\perp}Z\cdot \nabla \psi \right) \simeq \dv \left(\nabla \psi \Delta \psi\right)$ in terms of regularity and decay rates. So, the asymptotic behavior of $B$ is 
\[
B(t,x)=\left(\gamma \nabla^{\perp}\Gamma(t,x), \eta\Gamma(t,x)\right)+O(t^{-3/2}).
\]

\subsection{\bf Perturbation around harmonic functions} \label{sec:1.5}
Theorem \ref{GWP} is about the existence of a solution globally-in-time when the initial data is  small enough around zero. We now perturb (\ref{Hall equation in 2D}) around harmonic functions. Then, the newly generated terms are linear and so one may guess that a smallness condition on harmonic functions is also necessary. But, we emphasize that this is not the only case. As one can see from the statements of Theorem \ref{perturbation theorem 1} and Theorem \ref{perturbation theorem 2} or the proof of them in Section \ref{sec:4},  we can absorb these terms to the left-hand side of the desired bounds by multiplying by a large constant depending on  the harmonic function we choose.

\subsubsection{\bf Case 1}
Let $\overline{\psi}$ be a harmonic function such that $\left\|\nabla^{2}\overline{\psi}\right\|_{L^{\infty}}<\infty$. (For example, $\overline{\psi}(x,y)=x^{2}-y^{2}$ or $\overline{\psi}(x,y)=xy$.) Let $\psi=\rho+\overline{\psi}$. Then, we obtain the following equations of $(\rho, Z)$:
\begin{subequations} \label{new perturbation equation 1}
\begin{align}
&\rho_{t}-\Delta \rho=[\rho,Z]+ [\overline{\psi},Z], \\
& Z_{t} -\Delta Z=[\Delta \rho,\rho]+ [\Delta \rho, \overline{\psi}].  
\end{align}
\end{subequations}
Since the smallness condition is stated by combining $\overline{\psi}$ and the regularity of the various level of the initial conditions as shown just below, we take 
\begin{equation} \label{several norm 1}
\begin{split}
F_{1}=\left\|\nabla \rho\right\|^{2}_{L^{2}}+\left\|Z\right\|^{2}_{L^{2}}, \quad &F_{2}=\left\|\Delta \rho\right\|^{2}_{L^{2}}+\left\|\nabla Z\right\|^{2}_{L^{2}}\\
F_{3}=\left\|\nabla \Delta \rho\right\|^{2}_{L^{2}}+\left\|\Delta Z\right\|^{2}_{L^{2}}, \quad &F_{4}=\left\|\Delta^{2} \rho\right\|^{2}_{L^{2}}+\left\|\nabla \Delta Z\right\|^{2}_{L^{2}}
\end{split}
\end{equation}
and let $\epsilon_{2}=C_{1}^{2}F_{1}(0)+C_{1}F_{2}(0)+ F_{3}(0)$ and $C_{1}=k\left\|\nabla^{2}\overline{\psi}\right\|^{2}_{L^{\infty}}$ with $k$ fixed by (\ref{size k 1}).

\begin{theorem} \label{perturbation theorem 1}
There exists a constant $C$ such that if $C\epsilon_{2}<1$, there exists a unique solution $(\nabla \rho, Z)$ of (\ref{new perturbation equation 1}) satisfying  
\[
\begin{split}
&C_{1}^{2}F_{1}(t)+C_{1}F_{2}(t)+ F_{3}(t)+ (1-C\epsilon_{1})\int^{t}_{0}\left(C_{1}F_{3}(s)+F_{4}(s)\right)ds \\
&\leq C_{1}^{2}F_{1}(0)+C_{1}F_{2}(0)+ F_{3}(0)\quad \text{for all $t>0$.}
\end{split}
\]
\end{theorem}

\subsubsection{\bf Case 2}
Let $\overline{Z}$ be a harmonic function such that $\left\|\nabla \overline{Z}\right\|_{L^{\infty}}<\infty$.  (For example, $\overline{Z}(x,y)=ax+by$.) Let $Z=\omega+\overline{Z}$. Then, we obtain the following equations of $(\psi,\omega)$:
\begin{subequations} \label{new perturbation equation 2}
\begin{align}
&\psi_{t}-\Delta \psi=[\psi,\omega]+ [\psi, \overline{Z}], \\
& \omega_{t} -\Delta \omega=[\Delta \psi,\psi].  
\end{align}
\end{subequations}
For the same reason as Case 1, let  
\begin{equation} \label{several norm 2}
\begin{split}
K_{1}=\left\|\nabla \psi\right\|^{2}_{L^{2}}+\left\|\omega\right\|^{2}_{L^{2}}, \quad  &K_{2}=\left\|\Delta \psi\right\|^{2}_{L^{2}}+\left\|\nabla \omega\right\|^{2}_{L^{2}}, \\
K_{3}=\left\|\nabla \Delta \psi\right\|^{2}_{L^{2}}+\left\|\Delta \omega\right\|^{2}_{L^{2}}, \quad &K_{4}=\left\|\Delta^{2} \psi\right\|^{2}_{L^{2}}+\left\|\nabla \Delta \omega\right\|^{2}_{L^{2}}
\end{split}
\end{equation}
and let $\epsilon_{3}=C_{2}^{3}\left\|\psi_{0}\right\|^{2}_{L^{2}}+ C_{2}^{2}K_{1}(0)+ C_{2}K_{2}(0)+K_{3}(0)$  and $C_{2}=k\left\|\nabla^{2}\overline{Z}\right\|^{2}_{L^{\infty}}$ with $k$ fixed by (\ref{size k 2}). 

\begin{theorem} \label{perturbation theorem 2}
There exists a constant $C$ such that if $C\epsilon_{3}<1$, there exists a unique solution $(\nabla \psi, \omega)$ of (\ref{new perturbation equation 2}) satisfying 
\begin{equation*} \label{final bound per 2}
\begin{split}
&C_{2}^{3}\left\|\psi(t)\right\|^{2}_{L^{2}}+ C_{2}^{2}K_{1}(t)+C_{2}K_{2}(t)+K_{3}(t)+(1-C\epsilon_{2})\int^{t}_{0}\left(C_{2}^{2}K_{2}(s)+ C_{2}K_{3}(s)+K_{4}(s)\right)ds\\
&\leq C_{2}^{3}\left\|\psi(0)\right\|^{2}_{L^{2}}+ C_{2}^{2}K_{1}(0)+C_{2}K_{2}(0)+ K_{3}(0) \quad \text{for all $t>0$.}
\end{split}
\end{equation*}

\end{theorem}

\begin{remark} \upshape
\begin{enumerate}[]
\item (1) Compared to Theorem \ref{GWP}, we are not able to derive decay rates in Theorem \ref{perturbation theorem 1} and Theorem \ref{perturbation theorem 2} due to the terms having $\overline{\psi}$ and $\overline{Z}$ on the right-hand side of (\ref{perp 2}), (\ref{perp 3}), (\ref{new perp 2}), and (\ref{new perp 3}). 
\item (2) If we take $\psi=\rho+\overline{\psi}$ and $Z=\omega+\overline{Z}$, we obtain the following equations:
\begin{subequations} \label{both harmonic perturbation}
\begin{align}
&\rho_{t}-\Delta \rho=[\rho,\omega]+ [\rho, \overline{Z}]+[\overline{\psi},\omega]+ [\overline{\psi}, \overline{Z}], \\
& \omega_{t} -\Delta \omega=[\Delta \rho,\rho]+ [\Delta \rho, \overline{\psi}].
\end{align}
\end{subequations}
To obtain a result like Theorem \ref{perturbation theorem 1} and Theorem \ref{perturbation theorem 2}, we may start by choosing $\overline{\psi}$ and $\overline{Z}$ such that $[\overline{\psi}, \overline{Z}]=0$. One specific choice is $\overline{\psi}=\overline{Z}=ax+by$. In this case, $[\rho, \overline{Z}]=b\rho_{x}-a\rho_{y}$, $[\overline{\psi},\omega]=a\omega_{y}-b\omega_{x}$ and $[\Delta \rho, \overline{\psi}]=b\Delta \rho_{x}-a\Delta \rho_{y}$, and these will be cancelled out when we apply our method. So, (\ref{both harmonic perturbation}) is equivalent to (\ref{Hall equation in 2D}). But, we do not have results dealing with more general $\overline{\psi}$ and $\overline{Z}$ to handle  (\ref{both harmonic perturbation}) as (\ref{new perturbation equation 1}) or (\ref{new perturbation equation 2}).
\end{enumerate}
\end{remark}

\subsection{\bf Hall MHD}
After considering the Hall equations, we consider in this section the $2 \frac{1}{2}$ dimensional Hall MHD given by  (\ref{coupled two half}). Due to the presence of the fluid part, results similar to Theorem \ref{decay of psi}, Theorem \ref{Asymptotics},  Theorem \ref{perturbation theorem 1}, and Theorem \ref{perturbation theorem 2}  will not be presented in this paper. Instead, we begin with the existence and the decay rate of weak solutions of (\ref{coupled two half}) which are again the two dimensional version of \cite{Chae-Degond-Liu, Chae Schonbek}.

\begin{theorem}
Let $\left(\nabla\psi_{0}, Z_{0}, \nabla\phi_{0}, Z_{0}\right)\in L^{2}$. Then, there is a weak solution of (\ref{coupled two half}) satisfying 
\[
\begin{split}
&\left\|\nabla \psi(t)\right\|^{2}_{L^{2}}+\left\|Z(t)\right\|^{2}_{L^{2}}+\left\|\nabla \phi(t)\right\|^{2}_{L^{2}}+\left\|W(t)\right\|^{2}_{L^{2}}\\
&+2\int^{t}_{0}\left(\left\|\Delta \psi(s)\right\|^{2}_{L^{2}}+\left\|\nabla Z(s)\right\|^{2}_{L^{2}}+\left\|\Delta \phi(s)\right\|^{2}_{L^{2}}+\left\|\nabla W(s)\right\|^{2}_{L^{2}}\right)ds \\
&\leq \left\|\nabla \psi_{0}\right\|^{2}_{L^{2}}+\left\|Z_{0}\right\|^{2}_{L^{2}} +\left\|\nabla \phi_{0}\right\|^{2}_{L^{2}}+\left\|W_{0}\right\|^{2}_{L^{2}} 
\end{split}
\]
for all $t>0$. If $\left(\nabla\psi_{0}, Z_{0}, \nabla\phi_{0}, W_{0}\right)\in L^{2}\cap L^{1}$,  $\left(\nabla\psi, Z, \nabla\phi, W\right)$ decay in time as
\eqn\label{Weak decay Hall MHD}
\left\|\nabla \psi(t)\right\|_{L^{2}}+\left\|Z(t)\right\|_{L^{2}}+\left\|\nabla \phi(t)\right\|_{L^{2}}+\left\|W(t)\right\|_{L^{2}}\leq \frac{C_{0}}{\sqrt{1+t}},
\een
where $C_{0}$ depends on $\left\|\nabla \psi_{0}\right\|_{L^{2}\cap L^{1}}$, $\left\|Z_{0}\right\|_{L^{2}\cap L^{1}}$, $\left\|\nabla \phi_{0}\right\|_{L^{2}\cap L^{1}}$, and $\left\|W_{0}\right\|_{L^{2}\cap L^{1}}$.
\end{theorem}

We now proceed, as in Section \ref{sec:1.3}, towards the strong solutions of (\ref{coupled two half}). We first show 
the existence of unique local-in-time solutions with large initial data and we derive a blow-up criterion.  The function spaces that we introduce are similar to those used for (\ref{Hall equation in 2D})
\begin{equation}\label{energy norm 2}
\begin{split}
P(t)&=\left\|\nabla \psi(t)\right\|^{2}_{H^{2}}+\left\|Z(t)\right\|^{2}_{H^{2}}+\left\|\nabla \phi(t)\right\|^{2}_{H^{2}}+ \left\|W(t)\right\|^{2}_{H^{2}}, \\
Q(t)&=\left\|\Delta \psi(t)\right\|^{2}_{H^{2}}+\left\|\nabla Z(t)\right\|^{2}_{H^{2}}+\left\|\Delta \phi(t)\right\|^{2}_{H^{2}}+ \left\|\nabla W(t)\right\|^{2}_{H^{2}},\\
\mathcal{E}(t)&=P(t)+\int^{t}_{0}Q(s)ds.
\end{split}
\end{equation}
As for the Hall equations, constants that depend on $P(0)$ are not specifically specified each time when we state our results, and we will use $\mathcal{E}_{0}$ in common.  
 
\begin{theorem} \label{LWP Hall MHD}
Let $\left(\nabla\psi_{0}, Z_{0},\nabla \phi_{0}, W_{0}\right)\in H^{2}$. There exists $T^{\ast}=T(\mathcal{E}_{0})>0$ such that there exists a unique solution  of (\ref{coupled two half}) with $\mathcal{E}(T^{\ast})<\infty$. Moreover, the maximal existence time $T^{\ast}<\infty$ if and only if 
\eqn \label{blowup Hall MHD}
\lim_{T\nearrow T^{\ast}}\int^{T}_{0}\left\|\nabla Z(t)\right\|^{q}_{L^{p}}dt=\infty, \quad \frac{1}{p}+\frac{1}{q}=\frac{1}{2}, \quad 2\leq q<\infty.
\een
\end{theorem}

\begin{remark} \upshape
We emphasize that the blow-up criterion is specified only in terms of $Z$ even when the fluid part enters.
\end{remark}

In Section \ref{sec:7}, we will derive inequalities similar to (\ref{two inequalities Hall}) but the smallness condition is expressed more complicated by
\eqn \label{smallness condition 11}
\epsilon_{4}=\left\|\nabla \psi_{0}\right\|^{2}_{H^{1}}+\left\|Z_{0}\right\|^{2}_{H^{1}}+\left\|\nabla \phi_{0}\right\|^{2}_{H^{1}}+\left\|W_{0}\right\|^{2}_{H^{1}}, \quad C\epsilon_{4}<1.
\een
Compared to Theorem \ref{GWP}, we need to modify the smallness condition as (\ref{smallness condition 11}) because  (\ref{coupled two half}) does not have a scaling-invariant property. Suppose that $(u,B=0)$ and $(u=0,B)$ solves (\ref{Hall MHD}), respectively. Then, the same is true for rescaled functions: $u_{\lambda}(t,x)=\lambda u(\lambda^{2}t, \lambda x)$ and $B_{\lambda}(t,x)=B(\lambda^{2}t, \lambda x)$ accordingly. So, $u$ and $B$ have different scaling. Since (\ref{Hall MHD}) and so (\ref{coupled two half}) include both $u$ and $B$, the smallness condition is determined by a combination the scaling invariant quantities of $u$ and $B$.

\begin{theorem}\label{GWP Hall MHD}
Let $\left(\nabla\psi_{0}, Z_{0},\nabla \phi_{0}, W_{0}\right)\in H^{2}$ which satisfies (\ref{smallness condition 11}). Then, we can take $T=\infty$ in Theorem \ref{LWP Hall MHD}. $\left(\nabla\psi_{0}, Z_{0},\nabla \phi_{0}, W_{0}\right)\in L^{1}$ in addition, $(\Delta \psi, \nabla Z,\Delta \phi, \nabla W)$ decay in time as follows
\begin{equation} \label{GWP Decay Hall MHD}
\begin{split}
&\left\|\Delta \psi(t)\right\|_{L^{2}}+\left\|\nabla Z(t)\right\|_{L^{2}} +\left\|\Delta \psi(t)\right\|_{L^{2}}+\left\|\nabla Z(t)\right\|_{L^{2}}\leq \frac{\mathcal{E}_{0}}{1+t},\\
& \left\|\nabla \Delta \psi(t)\right\|_{L^{2}}+\left\|\Delta Z(t)\right\|_{L^{2}}+\left\|\nabla \Delta \phi(t)\right\|_{L^{2}}+\left\|\Delta W(t)\right\|_{L^{2}}\leq \frac{\mathcal{E}_{0}}{(1+t)^{3/2}}.
\end{split}
\end{equation}
\end{theorem}

The decay rate (\ref{GWP Decay Hall MHD}) can be easily derived by using (\ref{Weak decay Hall MHD}) and the argument in Section \ref{sec:3.3}, we will skip the process of proving (\ref{GWP Decay Hall MHD}).

\section{Preliminaries}

All  constants will be denoted by $C$ and we follow the convention that such constants can vary from expression to expression and even between two occurrences within the same expression. 

We here provide some inequalities in 2D:
\begin{subequations} \label{inequality 2}
\begin{align}
&\left\|f\right\|_{L^{p}}\leq C(p)\left\|f\right\|^{\frac{2}{p}}_{L^{2}}\left\|\nabla f\right\|^{1-\frac{2}{p}}_{L^{2}}, \quad 2<p<\infty\label{inequality 2 a}\\
& \left\|f\right\|_{L^{\infty}} \leq C\left\|f\right\|^{\frac{1}{2}}_{L^{2}} \left\|\Delta f\right\|^{\frac{1}{2}}_{L^{2}}\label{inequality 2 b}.
\end{align}
\end{subequations}
We also use the following inequalities which hold in any dimensions: 
\eqn \label{inequality 1}
\left\|\nabla f\right\|_{L^{2}} \leq  \left\|f\right\|^{\frac{1}{2}}_{L^{2}}\left\|\Delta f\right\|^{\frac{1}{2}}_{L^{2}}, \quad \left\|\nabla^{2}f\right\|_{L^{p}}\leq C(p)\left\|\Delta f\right\|_{L^{p}}, \quad 1<p<\infty.
\een

We now recall the commutator $[f,g]=\nabla f\cdot \nabla^{\perp}g=f_{x}g_{y}-f_{y}g_{x}$. Then, the commutator has the following properties:
\begin{subequations} \label{commutator}
\begin{align}
& [f,f]=0,  \quad [f,g]=-[g,f]\\
& \Delta [f,g]= [\Delta f,g] +[f,\Delta g] + 2[f_{x},g_{x}] +2[f_{y},g_{y}],\label{commutator c}\\
& \int f[f,g]=0, \label{commutator d}\\
& \int f[g,h]=\int g [h,f]. \label{commutator e}
\end{align}
\end{subequations}

We will use (\ref{inequality 2}) -- (\ref{commutator}) repeatedly when proving our results and we will not refer them every time when it is obvious to use them.

\subsection*{How to prove our results}
To prove our results, we may use  a fixed point argument. But, in principle, the calculations used to derive a priori estimates are easily applied to the fixed point argument. So, we only provide a priori estimates for the existence part and show the uniqueness. 	

\section{Proofs of Theorem \ref{LWP}, Theorem \ref{weak strong uniqueness}, and Theorem \ref{GWP}} \label{sec:3}
In this section we establish the local-in-time existence of unique strong solutions of (\ref{Hall equation in 2D}).  The analysis given in this section will apply to (\ref{coupled two half}) in Section \ref{sec:4} and Section \ref{sec:7} . Since the computations used to prove Theorem \ref{LWP} can be used to prove Theorem \ref{weak strong uniqueness}, we begin with Theorem \ref{LWP}. 

\subsection{Proof of Theorem \ref{LWP}}
We first recall (\ref{Hall equation in 2D}):
\begin{subequations}\label{equation psi Z}
\begin{align}
&\psi_{t}-\Delta \psi=[\psi,Z], \label{equation of psi a}\\
&Z_{t} -\Delta Z=[\Delta \psi,\psi].   \label{equation of Z a}
\end{align}
\end{subequations}

\subsubsection{\bf A priori estimates}\label{sec:3.1.1}
We multiply (\ref{equation of psi a}) by $-\Delta \psi$,  (\ref{equation of Z a}) by $Z$, and integrate over $\mathbb{R}^{2}$. By using (\ref{commutator e}), we have 
\eqn \label{L2 bound}
\frac{1}{2}\frac{d}{dt}\left(\left\|\nabla \psi\right\|^{2}_{L^{2}}+\left\|Z\right\|^{2}_{L^{2}}\right)+\left\|\Delta \psi\right\|^{2}_{L^{2}}+\left\|\nabla Z\right\|^{2}_{L^{2}}=-\int \Delta \psi[\psi,Z] +\int Z[\Delta \psi,\psi]=0.
\een

We next multiply (\ref{equation of psi a}) by $\Delta^{2} \psi$,  (\ref{equation of Z a}) by $-\Delta Z$ and integrate over $\mathbb{R}^{2}$. Then,  
\[
\frac{1}{2}\frac{d}{dt}\left(\left\|\Delta \psi\right\|^{2}_{L^{2}}+\left\|\nabla Z\right\|^{2}_{L^{2}}\right)+\left\|\nabla \Delta \psi\right\|^{2}_{L^{2}}+\left\|\Delta Z\right\|^{2}_{L^{2}} = \int \Delta^{2}  \psi [\psi,Z] -\int\Delta Z[\Delta \psi,\psi].
\] 
Since 
\begin{equation}\label{H2 bound 1}
\begin{split}
&\int \Delta^{2}  \psi [\psi,Z] -\int\Delta Z[\Delta \psi,\psi]=2\int \Delta \psi \left([\psi_{x},Z_{x}]+[\psi_{y},Z_{y}] \right)\leq C\left\|\nabla^{2}Z\right\|_{L^{2}} \left\|\nabla^{2}\psi\right\|^{2}_{L^{4}}\\
& \leq C \left\|\Delta Z\right\|_{L^{2}} \left\|\nabla^{2}\psi\right\|_{L^{2}}\left\|\nabla^{3}\psi\right\|_{L^{2}} \leq C\left\|\Delta \psi\right\|^{2}_{L^{2}}\left\|\nabla\Delta \psi\right\|^{2}_{L^{2}} +\frac{1}{2}\left\|\Delta Z\right\|^{2}_{L^{2}},
\end{split}
\end{equation}
we obtain
\eqn \label{H2 bound 2}
\frac{d}{dt}\left(\left\|\Delta \psi\right\|^{2}_{L^{2}}+\left\|\nabla Z\right\|^{2}_{L^{2}}\right)+\left\|\nabla \Delta \psi\right\|^{2}_{L^{2}}+\left\|\Delta Z\right\|^{2}_{L^{2}}\leq C\left\|\Delta \psi\right\|^{2}_{L^{2}} \left\|\nabla\Delta \psi\right\|^{2}_{L^{2}}.
\een

We finally multiply (\ref{equation of psi a}) by $-\Delta^{3} \psi$,  (\ref{equation of Z a}) by $\Delta^{2} Z$ and integrate over $\mathbb{R}^{2}$. By noticing, as (\ref{H2 bound 1}), the cancellation of the terms having the highest order derivative in the first equality below, we have  
\begin{equation*}
\begin{split}
&\frac{1}{2}\frac{d}{dt}\left(\left\|\nabla \Delta \psi\right\|^{2}_{L^{2}}+\left\|\Delta Z\right\|^{2}_{L^{2}}\right) +\left\|\Delta^{2} \psi\right\|^{2}_{L^{2}}+\left\|\nabla \Delta Z\right\|^{2}_{L^{2}} = -\int \Delta^{3}\psi [\psi,Z]+\int \Delta^{2} Z[\Delta \psi,\psi]\\
& = -\int \Delta^{2}\psi [\Delta \psi,Z] -2\int \Delta^{2}\psi \left([\psi_{x},Z_{x}]+[\psi_{y},Z_{y}]\right) - 2\int \Delta\psi \left([\psi_{x},\Delta Z_{x}]+[\psi_{y},\Delta Z_{y}]\right)\\
&=\text{(I)+(II)+(III)}.
\end{split}
\end{equation*}
We bound each term on the right-hand side as follows. By using the definition of the commutator,
\begin{equation} \label{first term}
\begin{split}
\text{(I)}&=-\int \left(\nabla^{\perp}Z\cdot \nabla \Delta \psi\right) \Delta^{2}\psi =\int \left(\nabla \nabla^{\perp}Z\cdot \nabla \Delta \psi\right)\cdot  \nabla \Delta\psi +\int \left(\nabla^{\perp}Z\cdot \nabla \nabla \Delta \psi\right) \nabla \Delta\psi \\
&=\int \left(\nabla \nabla^{\perp}Z\cdot \nabla \Delta \psi\right)\cdot  \nabla \Delta\psi -\frac{1}{2}\int \nabla \nabla^{\perp}Z \left|\nabla \Delta \psi\right|^{2} \leq C\int \left|\nabla^{2}Z\right| \left|\nabla^{3}\psi \right|^{2}\\
&\leq C \left\|\nabla^{2}Z\right\|_{L^{2}} \left\|\nabla^{3}\psi\right\|^{2}_{L^{4}}  \leq C\left\|\nabla^{2}Z\right\|_{L^{2}} \left\|\nabla^{3}\psi\right\|_{L^{2}} \left\|\nabla^{4}\psi\right\|_{L^{2}} \\
&\leq C\left\|\Delta Z\right\|^{2}_{L^{2}} \left\|\nabla\Delta \psi\right\|^{2}_{L^{2}}+\frac{1}{4}\left\|\Delta^{2}\psi\right\|^{2}_{L^{2}} .
\end{split}
\end{equation}
By moving one derivative in $\Delta^{2}\psi$ to $\left([\psi_{x},Z_{x}]+[\psi_{y},Z_{y}]\right)$ and by using (\ref{commutator e}), 
\begin{equation} \label{H3 bound 1}
\begin{split}
\text{(II)}+\text{(III)}\leq C\int \left|\nabla^{2}Z\right| \left|\nabla^{3}\psi \right|^{2}+C\int  \left|\nabla^{2}\psi \right| \left|\nabla^{3}\psi \right|\left|\nabla^{3}Z\right|
\end{split}
\end{equation}
with the second term estimated by 
\begin{equation*} 
\begin{split}
\int  \left|\nabla^{2}\psi \right| \left|\nabla^{3}\psi \right|\left|\nabla^{3}Z\right| &\leq \left\|\nabla^{2}\psi\right\|_{L^{4}} \left\|\nabla^{3}\psi\right\|_{L^{4}} \left\|\nabla^{3}Z\right\|_{L^{2}} \\
& \leq C \left\|\Delta \psi\right\|^{2}_{L^{2}} \left\|\nabla \Delta \psi\right\|^{4}_{L^{2}}+\frac{1}{4}\left\|\Delta^{2}\psi\right\|^{2}_{L^{2}}+ \frac{1}{2}\left\|\nabla\Delta Z\right\|^{2}_{L^{2}}.
\end{split}
\end{equation*}
With these estimates, we have  
\begin{equation}\label{H3 bound 2}
\begin{split}
&\frac{d}{dt}\left(\left\|\nabla \Delta \psi\right\|^{2}_{L^{2}}+\left\|\Delta Z\right\|^{2}_{L^{2}}\right) +\left\|\Delta^{2} \psi\right\|^{2}_{L^{2}}+\left\|\nabla \Delta Z\right\|^{2}_{L^{2}} \\
&\leq  C\left\|\Delta Z\right\|^{2}_{L^{2}} \left\|\nabla\Delta \psi\right\|^{2}_{L^{2}}+C \left\|\Delta \psi\right\|^{2}_{L^{2}} \left\|\nabla \Delta \psi\right\|^{4}_{L^{2}}.
\end{split}
\end{equation}

By (\ref{L2 bound}), (\ref{H2 bound 2}), and (\ref{H3 bound 2}), we derive the following inequality:  
\eqn \label{full H3 ddd}
\frac{d}{dt}(1+M)+N\leq CM^{2}+C M^{3}\leq C(1+M)^{3}, 
\een
where $M$ and $N$ are defined in (\ref{energy norm}). From this, we deduce
\eqn \label{full H3 dddd}
M(t)\leq \sqrt{\frac{(1+M(0))^{2}}{1-2Ct(1+M(0))^{2}}}-1 \quad \text{for all} \ t\leq T^{\ast}<\frac{1}{2C (1+M(0))^{2}}.
\een
Integrating (\ref{full H3 ddd}) and using (\ref{full H3 dddd}), we finally derive $\mathcal{E}(T^{\ast})<\infty$.

\subsubsection{\bf Uniqueness} \label{sec:3.1.2}
Suppose there are two solutions $(\psi_{1}, Z_{1})$ and $(\psi_{2}, Z_{2})$. Let $\psi=\psi_{1}-\psi_{2}$ and $Z=Z_{1}-Z_{2}$. By subtracting the equations for $(\psi_{1}, Z_{1})$ and $(\psi_{2}, Z_{2})$, we have
\begin{subequations} \label{difference}
\begin{align}
&\psi_{t}-\Delta \psi=[\psi_{1}, Z]+[\psi,Z_{2}], \\
& Z_{t}-\Delta Z=[\Delta \psi, \psi_{1}]+[\Delta \psi_{2}, \psi].
\end{align}
\end{subequations}
From this, we see that   
\begin{equation}  \label{difference uniqueness}
\begin{split}
&\frac{1}{2}\frac{d}{dt} \left(\left\|\nabla \psi\right\|^{2}_{L^{2}}+\left\|Z\right\|^{2}_{L^{2}}\right) +\left\|\Delta \psi\right\|^{2}_{L^{2}}+ \left\|\nabla Z\right\|^{2}_{L^{2}}\\
& =-\int \Delta \psi[\psi_{1},Z] - \int \Delta \psi[\psi,Z_{2}] +\int Z[\Delta \psi,\psi_{1}] +\int Z[\Delta \psi_{2},\psi]\\
&= - \int \Delta \psi[\psi,Z_{2}]  +\int Z[\Delta \psi_{2},\psi] \\
&\leq C\left(\left\|\nabla Z_{2}\right\|^{2}_{L^{2}} \left\|\Delta Z_{2}\right\|^{2}_{L^{2}} +\left\|\Delta \psi_{2}\right\|^{2}_{L^{2}} \left\|\nabla \Delta \psi_{2}\right\|^{2}_{L^{2}} \right)\left\|\nabla \psi\right\|^{2}_{L^{2}} +\left\|\Delta \psi\right\|^{2}_{L^{2}}+ \left\|\nabla Z\right\|^{2}_{L^{2}}
\end{split}
\end{equation}
which gives 
\[
\frac{d}{dt} \left(\left\|\nabla \psi\right\|^{2}_{L^{2}}+\left\|Z\right\|^{2}_{L^{2}}\right) \leq C\left(\left\|\nabla Z_{2}\right\|^{2}_{L^{2}} \left\|\Delta Z_{2}\right\|^{2}_{L^{2}} +\left\|\Delta \psi_{2}\right\|^{2}_{L^{2}} \left\|\nabla \Delta \psi_{2}\right\|^{2}_{L^{2}} \right)\left(\left\|\nabla \psi\right\|^{2}_{L^{2}}+\left\|Z\right\|^{2}_{L^{2}}\right). 
\]
Since $\left\|\nabla Z_{2}\right\|^{2}_{L^{2}} \left\|\Delta Z_{2}\right\|^{2}_{L^{2}} +\left\|\Delta \psi_{2}\right\|^{2}_{L^{2}} \left\|\nabla \Delta \psi_{2}\right\|^{2}_{L^{2}}$  is integrable on $[0,T^{\ast})$, the uniqueness follows using Gronwall's lemma.

\subsubsection{\bf Blow-up criterion} \label{sec:3.1.3}
To  obtain (\ref{blowup}), we first bound the right-hand side of (\ref{H2 bound 1}) by
\[
\left|\int \Delta \psi \left([\psi_{x},Z_{x}]+[\psi_{y},Z_{y}] \right)\right|\leq C\left\|\nabla Z\right\|_{L^{p}} \left\|\nabla^{2}\psi\right\|_{L^{q}} \left\|\nabla^{3}\psi\right\|_{L^{2}}, \quad \frac{1}{p}+\frac{1}{q}=\frac{1}{2}.
\]
By (\ref{inequality 2 a}), we have
\begin{equation} \label{blowup first}
\begin{split}
\left\|\nabla Z\right\|_{L^{p}} \left\|\Delta\psi\right\|_{L^{q}} \left\|\nabla^{3}\psi\right\|_{L^{2}}&\leq C\left\|\nabla Z\right\|_{L^{p}} \left\|\nabla^{2}\psi\right\|^{\frac{2}{q}}_{L^{2}} \left\|\nabla\Delta\psi\right\|^{2-\frac{2}{q}}_{L^{2}}\\
&\leq C\left\|\nabla Z\right\|^{q}_{L^{p}} \left\|\Delta\psi\right\|^{2}_{L^{2}} +\left\|\nabla\Delta\psi\right\|^{2}_{L^{2}}. 
\end{split}
\end{equation}
So, we can rewrite (\ref{H2 bound 1}) as
\[ \label{blowup 1}
\frac{d}{dt}\left(\left\|\Delta \psi\right\|^{2}_{L^{2}}+\left\|\nabla Z\right\|^{2}_{L^{2}}\right)+\left\|\nabla \Delta \psi\right\|^{2}_{L^{2}}+\left\|\Delta Z\right\|^{2}_{L^{2}}\leq C\left\|\nabla Z\right\|^{q}_{L^{p}}\left(\left\|\Delta \psi\right\|^{2}_{L^{2}}+\left\|\nabla Z\right\|^{2}_{L^{2}}\right). 
\]
Integrating this in time by using Gronwall's inequality, we have
\begin{equation} \label{blowup 2}
\begin{split}
&\left\|\Delta \psi(t)\right\|^{2}_{L^{2}}+\left\|\nabla Z(t)\right\|^{2}_{L^{2}}+\int^{t}_{0} \left(\left\|\nabla \Delta \psi(s)\right\|^{2}_{L^{2}}+\left\|\Delta Z(s)\right\|^{2}_{L^{2}}\right)ds\\
&\leq \mathcal{E}_{0}\exp\left[C\int^{t}_{0}\left\|\nabla Z(s)\right\|^{q}_{L^{p}}ds\right].
\end{split}
\end{equation}
We then integrate (\ref{H3 bound 2}) in time, without including $\left\|\Delta^{2} \psi\right\|^{2}_{L^{2}}+\left\|\nabla \Delta Z\right\|^{2}_{L^{2}}$, to obtain 
\[
\left\|\nabla \Delta \psi(t)\right\|^{2}_{L^{2}}+\left\|\Delta Z(t)\right\|^{2}_{L^{2}}\leq \mathcal{E}_{0}\exp\left[C\int^{t}_{0}\left(\left\|\nabla \Delta \psi\right\|^{2}_{L^{2}}+\left\|\Delta \psi\right\|^{2}_{L^{2}} \left\|\nabla \Delta \psi\right\|^{2}_{L^{2}}\right)ds\right].
\]
By using (\ref{blowup 2}) for the second term in the integrand of the right-hand side, we obtain 
\eqn \label{Blowup Hall equation}
\left\|\nabla\Delta \psi(t)\right\|^{2}_{L^{2}}+\left\|\Delta Z(t)\right\|^{2}_{L^{2}}\leq \mathcal{E}_{0}\exp \exp\left[CB(t)+CB^{2}(t)\right], \quad B(t)=\int^{t}_{0}\left\|\nabla Z(s)\right\|^{q}_{L^{p}}ds.
\een
This completes the proof of Theorem \ref{LWP}.

\subsection{Proof of Theorem \ref{weak strong uniqueness}}
The proof of Theorem \ref{weak strong uniqueness} is very similar to the one in Section \ref{sec:3.1.2} and Section \ref{sec:3.1.3}. Let $B_{1}=\left(\nabla^{\perp}\psi_{1}, Z_{1}\right)$ and $B_{2}=\left(\nabla^{\perp}\psi_{2}, Z_{2}\right)$ be the two weak solutions of (\ref{Hall equation in 2D}). Let $\psi=\psi_{1}-\psi_{2}$ and $Z=Z_{1}-Z_{2}$. By (\ref{difference uniqueness}) and (\ref{commutator e}), and by using  (\ref{blowup first}), we have 
\begin{equation*}
\begin{split}
&\frac{1}{2}\frac{d}{dt} \left(\left\|\nabla \psi\right\|^{2}_{L^{2}}+\left\|Z\right\|^{2}_{L^{2}}\right) +\left\|\Delta \psi\right\|^{2}_{L^{2}}+ \left\|\nabla Z\right\|^{2}_{L^{2}}= - \int \Delta \psi[\psi,Z_{2}]  +\int \Delta \psi_{2}[\psi,Z]\\
&\leq C\left\|\nabla Z_{2}\right\|_{L^{p}} \left\|\nabla\psi\right\|_{L^{q}} \left\|\Delta\psi\right\|_{L^{2}}+C\left\|\Delta \psi_{2}\right\|_{L^{p}} \left\|\nabla\psi\right\|_{L^{q}} \left\|\nabla Z\right\|_{L^{2}}\\
& \leq C\left(\left\|\nabla Z_{2}\right\|^{q}_{L^{p}}+\left\|\Delta \psi_{2}\right\|^{q}_{L^{p}} \right)\left(\left\|\nabla \psi\right\|^{2}_{L^{2}}+\left\|Z\right\|^{2}_{L^{2}}\right)+\frac{1}{2}\left\|\Delta \psi\right\|^{2}_{L^{2}}+ \frac{1}{2}\left\|\nabla Z\right\|^{2}_{L^{2}}
\end{split}
\end{equation*}
and so we obtain 
\[
\frac{d}{dt} \left(\left\|\nabla \psi\right\|^{2}_{L^{2}}+\left\|Z\right\|^{2}_{L^{2}}\right) +\left\|\Delta \psi\right\|^{2}_{L^{2}}+ \left\|\nabla Z\right\|^{2}_{L^{2}}\leq C\left(\left\|\nabla Z_{2}\right\|^{q}_{L^{p}}+\left\|\Delta \psi_{2}\right\|^{q}_{L^{p}} \right)\left(\left\|\nabla \psi\right\|^{2}_{L^{2}}+\left\|Z\right\|^{2}_{L^{2}}\right).
\]
By Gronwall inequality, we complete the proof of Theorem \ref{weak strong uniqueness}.

\subsection{Proof of Theorem \ref{GWP}}\label{sec:3.3}

\subsubsection{\bf A priori estimates}
We now show that the strong solutions provided by Theorem \ref{LWP} are in fact defined for all $t>0$ under the smallness condition (\ref{smallness condition 1}). 
We first rewrite  (\ref{H2 bound 2}) as 
\[
\frac{d}{dt}\left(\left\|\Delta \psi\right\|^{2}_{L^{2}}+\left\|\nabla Z\right\|^{2}_{L^{2}}\right)+\left\|\nabla \Delta \psi\right\|^{2}_{L^{2}}+\left\|\Delta Z\right\|^{2}_{L^{2}}\leq CS(t)\left(\left\|\nabla \Delta \psi\right\|^{2}_{L^{2}}+\left\|\Delta Z\right\|^{2}_{L^{2}}\right),
\]
where $S(t)=\left\|\Delta \psi\right\|^{2}_{L^{2}}+\left\|\nabla Z\right\|^{2}_{L^{2}}$. Let $\epsilon_{1}=\left\|\Delta \psi_{0}\right\|^{2}_{L^{2}}+ \left\|\nabla Z_{0}\right\|^{2}_{L^{2}}$. If $C\epsilon_{1}<1$, we have
\eqn \label{Hall decay 1}
\left\|\Delta \psi(t)\right\|^{2}_{L^{2}}+\left\|\nabla Z(t)\right\|^{2}_{L^{2}}+(1-C\epsilon_{1})\int^{t}_{0}\left(\left\|\nabla \Delta \psi(s)\right\|^{2}_{L^{2}}+\left\|\Delta Z(s)\right\|^{2}_{L^{2}}\right)ds\leq \epsilon_{1}
\een
for all $t>0$. We next  proceed to bound  (\ref{H3 bound 2}) by estimating the two terms on the right-hand side of (\ref{H3 bound 1}) in a different way. From the last expression of (\ref{first term}), we obtain 
\begin{equation} \label{third term}
\begin{split}
&C\left\|\nabla^{2}Z\right\|_{L^{2}} \left\|\nabla^{3}\psi\right\|^{2}_{L^{4}}\leq C\left\|\Delta Z\right\|^{2}_{L^{2}} \left\|\nabla\Delta \psi\right\|^{2}_{L^{2}}+\frac{1}{2}\left\|\Delta^{2}\psi\right\|^{2}_{L^{2}} \\
& \leq C\left\|\nabla Z\right\|^{2}_{L^{2}} \left\|\nabla \Delta Z\right\|^{2}_{L^{2}} +C\left\|\Delta \psi\right\|^{2}_{L^{2}} \left\|\Delta^{2} \psi\right\|^{2}_{L^{2}}+\frac{1}{2}\left\|\Delta^{2}\psi\right\|^{2}_{L^{2}}.
\end{split}
\end{equation}
By (\ref{inequality 2 b}), we bound the second term on the right-hand side of (\ref{H3 bound 1}) as 
\begin{equation}\label{fourth term}
\begin{split}
&\left\|\nabla^{2}\psi\right\|_{L^{\infty}} \left\|\nabla^{3}\psi\right\|_{L^{2}} \left\|\nabla^{3}Z\right\|_{L^{2}} \leq C\left\|\Delta\psi\right\|^{\frac{1}{2}}_{L^{2}} \left\|\Delta^{2}\psi\right\|^{\frac{1}{2}}_{L^{2}} \left\|\nabla^{3}\psi\right\|_{L^{2}}  \left\|\nabla^{3}Z\right\|_{L^{2}}\\
& \leq C\left\|\Delta\psi\right\|_{L^{2}} \left\|\Delta^{2}\psi\right\|_{L^{2}}  \left\|\nabla^{3}Z\right\|_{L^{2}} \leq C\left\|\Delta\psi\right\|^{2}_{L^{2}} \left\|\Delta^{2}\psi\right\|^{2}_{L^{2}}+ \frac{1}{2}\left\|\nabla^{3}Z\right\|^{2}_{L^{2}}.
\end{split}
\end{equation}
So, (\ref{H3 bound 2}) is replaced with
\eqn  \label{H3 bound 3}
\frac{d}{dt}\left(\left\|\nabla \Delta \psi\right\|^{2}_{L^{2}}+\left\|\Delta Z\right\|^{2}_{L^{2}}\right) +\left\|\Delta^{2} \psi\right\|^{2}_{L^{2}}+\left\|\nabla \Delta Z\right\|^{2}_{L^{2}} \leq  CS(t)\left(\left\|\Delta^{2} \psi\right\|^{2}_{L^{2}}+\left\|\nabla \Delta Z\right\|^{2}_{L^{2}} \right).
\een
Then (\ref{Hall decay 1}) gives 
\begin{equation*} \label{Hall decay 2}
\begin{split}
&\left\|\nabla\Delta \psi(t)\right\|^{2}_{L^{2}}+\left\|\Delta Z(t)\right\|^{2}_{L^{2}}+(1-C\epsilon_{1})\int^{t}_{0}\left(\left\|\Delta^{2} \psi(s)\right\|^{2}_{L^{2}}+\left\|\nabla \Delta Z(s)\right\|^{2}_{L^{2}}\right)ds\\
&\leq \left\|\nabla\Delta \psi_{0}\right\|^{2}_{L^{2}}+\left\|\Delta Z_{0}\right\|^{2}_{L^{2}}
\end{split}
\end{equation*}
for all $t>0$. This completes the first part of Theorem \ref{GWP}.

\subsubsection{\bf Decay rates} \label{sec:3.3.2}
To conclude this Section and the proof of Theorem \ref{GWP}, we now prove the decay rates (\ref{GWP Decay 2}) in Theorem \ref{GWP}. We first write (\ref{H2 bound 2}) as 
\eqn \label{H2 bound 3}
\frac{d}{dt} \left(\left\|\Delta \psi\right\|^{2}_{L^{2}}+\left\|\nabla Z\right\|^{2}_{L^{2}} \right) +(1-C\epsilon_{1})\left\|\nabla\Delta \psi\right\|^{2}_{L^{2}} +(1-C\epsilon_{1})\left\|\Delta Z\right\|^{2}_{L^{2}} \leq 0.
\een
By  (\ref{inequality 1}) and (\ref{Weak decay}), we have  
\begin{equation*}
\begin{split}
&\left\|\Delta \psi\right\|^{4}_{L^{2}}\leq  C\left\|\nabla \psi\right\|^{2}_{L^{2}}\left\|\nabla \Delta \psi\right\|^{2}_{L^{2}} \leq \frac{\mathcal{E}_{0}}{1+t} \left\|\nabla \Delta \psi\right\|^{2}_{L^{2}},\\
& \left\|\nabla Z\right\|^{4}_{L^{2}}\leq C \left\|Z\right\|^{2}_{L^{2}}\left\|\Delta Z\right\|^{2}_{L^{2}}\leq \frac{\mathcal{E}_{0}}{1+t}\left\|\Delta Z\right\|^{2}_{L^{2}}.
\end{split}
\end{equation*}
Then, (\ref{H2 bound 3}) becomes 
\[ 
\frac{d}{dt} \left(\left\|\Delta \psi\right\|^{2}_{L^{2}}+\left\|\nabla Z\right\|^{2}_{L^{2}} \right) +\mathcal{E}_{0}(1+t)\left(\left\|\Delta \psi\right\|^{2}_{L^{2}}+\left\|\nabla Z\right\|^{2}_{L^{2}} \right)^{2} \leq 0.
\]
By solving this ODE, we derive the following inequality for $t>0$:
\[ \label{decay 1}
\left\|\Delta \psi(t)\right\|^{2}_{L^{2}}+\left\|\nabla Z(t)\right\|^{2}_{L^{2}}\leq \frac{2\left\|\Delta \psi_{0}\right\|^{2}_{L^{2}}+ 2\left\|\nabla Z_{0}\right\|^{2}_{L^{2}}}{2+\mathcal{E}_{0}\big(\left\|\Delta \psi_{0}\right\|^{2}_{L^{2}}+ \left\|\nabla Z_{0}\right\|^{2}_{L^{2}}\big)(1+t)^{2}}.
\]

We next  write (\ref{H3 bound 3}) as 
\[
\frac{d}{dt} \left(\left\|\nabla \Delta \psi\right\|^{2}_{L^{2}}+\left\|\Delta Z\right\|^{2}_{L^{2}} \right) +\left(1-C\epsilon_{1}\right)\left\|\Delta^{2} \psi\right\|^{2}_{L^{2}} +\left(1-C\epsilon_{1}\right)\left\|\nabla \Delta Z\right\|^{2}_{L^{2}} \leq 0.
\]
By  (\ref{inequality 1}) and (\ref{Weak decay}), we have  
\begin{equation*}
\begin{split}
&\left\|\nabla \Delta \psi\right\|^{3}_{L^{2}} \leq C\left\|\nabla \psi\right\|_{L^{2}}\left\|\Delta^{2} \psi\right\|^{2}_{L^{2}}\leq \frac{\mathcal{E}_{0}}{\sqrt{1+t}}\left\|\Delta^{2} \psi\right\|^{2}_{L^{2}},\\
& \left\|\Delta Z\right\|^{3}_{L^{2}}\leq C \left\|Z\right\|_{L^{2}}\left\|\nabla \Delta Z\right\|^{2}_{L^{2}} \leq \frac{\mathcal{E}_{0}}{\sqrt{1+t}}\left\|\nabla \Delta Z\right\|^{2}_{L^{2}}.
\end{split}
\end{equation*}
So, we obtain  
\[ 
\frac{d}{dt} \left(\left\|\nabla \Delta \psi\right\|^{2}_{L^{2}}+\left\|\Delta Z\right\|^{2}_{L^{2}} \right) +\mathcal{E}_{0}\sqrt{1+t}\left(\left\|\nabla \Delta \psi\right\|^{2}_{L^{2}}+\left\|\Delta Z\right\|^{2}_{L^{2}} \right)^{\frac{3}{2}} \leq 0.
\]
From this, we derive the following inequality:
\[ \label{decay 2}
\left\|\nabla \Delta \psi(t)\right\|^{2}_{L^{2}}+\left\|\Delta Z(t)\right\|^{2}_{L^{2}}\leq \frac{36\left\|\nabla \Delta \psi_{0}\right\|^{2}_{L^{2}}+36\left\|\Delta Z_{0}\right\|^{2}_{L^{2}}}{\big(6+\mathcal{E}_{0}\sqrt{\left\|\nabla \Delta \psi_{0}\right\|^{2}_{L^{2}}+\left\|\Delta Z_{0}\right\|^{2}_{L^{2}}}(1+)^{\frac{3}{2}}\big)^{2}}
\] 
and thus concluding the proof of Theorem \ref{GWP}.

\section{Proof of Theorem \ref{decay of psi}}
In this section, we want to improve the decay rate of $\psi$ by using Theorem \ref{GWP}. We first recall the equation of $\psi$:
\eqn \label{eq of psi only}
\psi_{t}+\nabla^{\perp}Z\cdot \nabla \psi-\Delta \psi=0
\een
which is a dissipative transport equation with  a fast decaying coefficient $\nabla^{\perp}Z$. We begin with the $L^{1}$ bound of $\psi$:
\eqn \label{L1 bound}
\left\|\psi(t)\right\|_{L^1}\leq \left\|\psi_{0}\right\|_{L^1}.
\een
By applying Fourier splitting method in \cite{Chae Schonbek}, we also obtain the $L^{2}$ bound:
\eqn \label{L2 bound psi}
\left\|\psi(t)\right\|_{L^2}\leq \frac{\mathcal{E}_{0}}{\sqrt{1+t}}.
\een

We now test  $\displaystyle \frac{\psi}{t}$ to (\ref{eq of psi only}). Then, we obtain 
\begin{equation}\label{psi t}
\frac{1}{2}\frac{d}{dt}\frac{\left\|\psi\right\|^{2}_{L^{2}}}{t}+\frac{\left\|\psi\right\|^{2}_{L^{2}}}{2t^2}+\frac{\left\|\nabla \psi\right\|^{2}_{L^{2}}}{t}=0.
\end{equation}
By (\ref{inequality 2 a}) and (\ref{L1 bound}), we have 
\[
\left\|\psi\right\|^{4}_{L^{2}}\leq C \left\|\psi\right\|^{2}_{L^{1}}\left\|\nabla \psi\right\|^{2}_{L^{2}}\leq \mathcal{E}_{0} \left\|\nabla \psi\right\|^{2}_{L^{2}}
\]
and so (\ref{psi t}) can be replaced with 
\eqn \label{psi t 2}
\frac{1}{2}\frac{d}{dt}\frac{\left\|\psi\right\|^{2}_{L^{2}}}{t}+\frac{t}{\mathcal{E}_{0}}\Big(\frac{\left\|\psi\right\|^{2}_{L^{2}}}{t}\Big)^{2}\leq 0.
\een
On the other hand, testing $-\Delta \psi$ to (\ref{eq of psi only}), we have
\[
\frac{1}{2}\frac{d}{dt}\left\|\nabla \psi\right\|^{2}_{L^{2}}+\left\|\Delta \psi \right\|^{2}_{L^{2}}\leq \left\|\nabla Z\right\|_{L^4}\left\|\nabla \psi\right\|_{L^4}\left\|\Delta \psi\right\|_{L^{2}}\leq C\left\|\nabla Z\right\|^{4}_{L^4}\left\|\nabla \psi\right\|^{2}_{L^2}+\frac{1}{2}\left\|\Delta \psi\right\|^{2}_{L^{2}} 
\]
and so we obtain 
\eqn   \label{nabla psi t}
\frac{1}{2}\frac{d}{dt}\left\|\nabla \psi\right\|^{2}_{L^{2}}+\frac{1}{2}\left\|\Delta \psi \right\|^{2}_{L^{2}}\leq C\left\|\nabla Z\right\|^{4}_{L^4}\left\|\nabla \psi\right\|^{2}_{L^2}.
\een
By (\ref{inequality 2 a}) and (\ref{GWP Decay 2}),
\[
\left\|\nabla Z(t)\right\|^{4}_{L^4}\leq C \left\|\nabla Z(t)\right\|^{2}_{L^2}\left\|\Delta Z(t)\right\|^{2}_{L^2}\leq \frac{\mathcal{E}_{0}}{(1+t)^{5}}.
\]
By taking $t$ sufficiently large, which is expressed by $t>t_{0}$ for the rest of the proof of Theorem \ref{decay of psi},  and by combining (\ref{psi t 2}) and (\ref{nabla psi t}), we obtain
\[
\frac{1}{2}\frac{d}{dt}\Big(\frac{\left\|\psi\right\|^{2}_{L^{2}}}{t}+ \left\|\nabla \psi\right\|^{2}_{L^{2}}\Big)+\frac{1}{t}\left\|\nabla \psi\right\|^{2}_{L^{2}}+\left\|\Delta \psi \right\|^{2}_{L^{2}} \leq  0.
\]
By (\ref{inequality 1}) and (\ref{L2 bound psi}), we have 
\[
\left\|\nabla \psi\right\|^{4}_{L^{2}}\leq C \left\|\psi\right\|^{2}_{L^{2}}\left\|\Delta \psi\right\|^{2}_{L^{2}}\leq \frac{\mathcal{E}_{0}}{1+t} \left\|\Delta \psi\right\|^{2}_{L^{2}}\leq \frac{\mathcal{E}_{0}}{t} \left\|\Delta \psi\right\|^{2}_{L^{2}}
\]
when $t>t_{0}$. So, we derive the following inequality:
\eqn \label{t inequality}
\frac{1}{2}\frac{d}{dt}\Big(\frac{\left\|\psi\right\|^{2}_{L^{2}}}{t}+ \left\|\nabla \psi\right\|^{2}_{L^{2}}\Big)+\frac{t}{\mathcal{E}_{0}} \Big(\frac{\left\|\psi\right\|^{2}_{L^{2}}}{t}+ \left\|\nabla \psi\right\|^{2}_{L^{2}}\Big)^{2}\leq  0.
\een
We now solve this ODE to find
\[
\frac{\left\|\psi(t)\right\|^{2}_{L^{2}}}{t}+ \left\|\nabla \psi(t)\right\|^{2}_{L^{2}}\leq \frac{\mathcal{E}_{0}\Big(\frac{\left\|\psi(t_{0})\right\|^{2}_{L^{2}}}{t_{0}}+ \left\|\nabla \psi(t_{0})\right\|^{2}_{L^{2}}\Big)}{\mathcal{E}_{0} + \Big(\frac{\left\|\psi(t_{0})\right\|^{2}_{L^{2}}}{t_{0}}+ \left\|\nabla \psi(t_{0})\right\|^{2}_{L^{2}}\Big)(t^{2}-t^{2}_{0})}.
\]
Since $ \left\|\psi(t)\right\|_{H^{1}}\leq  \left\|\psi_{0}\right\|_{H^{1}}$ for all $t>0$ by Theorem \ref{GWP}, we obtain 
\eqn \label{decay 3}
\left\|\nabla \psi(t)\right\|_{L^{2}}\leq \frac{\mathcal{E}_{0}}{1+t}.
\een

\vspace{1ex}

By modifying (\ref{t inequality}) with the extra $t$-factor in the denominator, we have
\eqn \label{t inequality 2}
\frac{1}{2}\frac{d}{dt}\Big(\frac{\left\|\psi\right\|^{2}_{L^{2}}}{t^{2}}+ \frac{\left\|\nabla \psi\right\|^{2}_{L^{2}}}{t}\Big)+\frac{t^{2}}{\mathcal{E}_{0}} \Big(\frac{\left\|\psi\right\|^{2}_{L^{2}}}{t^{2}}+ \frac{\left\|\nabla \psi\right\|^{2}_{L^{2}}}{t}\Big)^{2}\leq  0.
\een 
We now test $\Delta^{2}\psi$ to (\ref{eq of psi only}). Then, we have
\begin{equation*}
\begin{split}
\frac{1}{2}\frac{d}{dt}\left\|\Delta \psi\right\|^{2}_{L^{2}}+\left\|\nabla \Delta \psi\right\|^{2}_{L^{2}}&=\int \Delta \psi \Delta [\psi,Z]=\int \Delta \psi [\psi,\Delta Z]+2\int \Delta \psi \left([\psi_{x}, Z_{x}]+[\psi_{y, Z_{y}}]\right)\\
&=\text{(I)+(II)}.
\end{split}
\end{equation*}
We first bound $\text{(I)}$ as follows
\begin{equation*}
\begin{split}
\text{(I)}&\leq \left\|\Delta Z\right\|_{L^{2}}\left\|\nabla \psi\right\|_{L^{\infty}}\left\|\nabla \Delta \psi\right\|_{L^{2}}\leq \left\|\Delta Z\right\|^{2}_{L^{2}}\left\|\nabla \psi\right\|^{2}_{L^{\infty}}+\frac{1}{8}\left\|\nabla \Delta \psi\right\|^{2}_{L^{2}}\\
&\leq C\left\|\Delta Z\right\|^{4}_{L^{2}}\left\|\nabla \psi\right\|^{2}_{L^{2}}+\frac{1}{4}\left\|\nabla \Delta \psi\right\|^{2}_{L^{2}}\leq \frac{\mathcal{E}_{0}}{(1+t)^{4}}\left\|\nabla \psi\right\|^{2}_{L^{2}} +\frac{1}{4}\left\|\nabla \Delta \psi\right\|^{2}_{L^{2}}.
\end{split}
\end{equation*}
$\text{(II)}$ is bounded by
\[
\begin{split}
\text{(II)}&\leq C \left\|\Delta Z\right\|_{L^{2}}\left\|\Delta \psi\right\|^{2}_{L^{4}}\leq C\left\|\Delta Z\right\|^{2}_{L^{2}}\left\|\Delta \psi\right\|^{2}_{L^{2}}+\frac{1}{4}\left\|\nabla \Delta \psi\right\|^{2}_{L^{2}}\\
&\leq \frac{\mathcal{E}_{0}}{(1+t)^{2}}\left\|\Delta \psi\right\|^{2}_{L^{2}} +\frac{1}{4}\left\|\nabla \Delta \psi\right\|^{2}_{L^{2}}.
\end{split}
\]
So, we obtain 
\eqn \label{Delta psi}
\frac{d}{dt}\left\|\Delta \psi\right\|^{2}_{L^{2}}+\left\|\nabla \Delta \psi\right\|^{2}_{L^{2}}\leq \frac{\mathcal{E}_{0}}{t^{4}}\left\|\nabla \psi\right\|^{2}_{L^{2}}+\frac{\mathcal{E}_{0}}{t^{2}}\left\|\Delta \psi\right\|^{2}_{L^{2}}
\een
when $t>t_{0}$. Then, (\ref{t inequality 2}) and (\ref{Delta psi}) give  
\[
\frac{d}{dt}\Big(\frac{\left\|\psi\right\|^{2}_{L^{2}}}{t^{2}}+ \frac{\left\|\nabla \psi\right\|^{2}_{L^{2}}}{t}\Big)+\mathcal{E}_{0}t^{2}\Big(\frac{\left\|\psi\right\|^{2}_{L^{2}}}{t^{2}}+ \frac{\left\|\nabla \psi\right\|^{2}_{L^{2}}}{t}\Big)^{2}+\left\|\nabla \Delta \psi\right\|^{2}_{L^{2}}\leq 0.
\]
By (\ref{inequality 1}) and  (\ref{decay 3}), we have
\[
\left\|\Delta \psi\right\|^{4}_{L^{2}}\leq C\left\|\nabla \psi\right\|^{2}_{L^{2}} \left\|\nabla \Delta \psi\right\|^{2}_{L^{2}}\leq \frac{\mathcal{E}_{0}}{(1+t)^{2}}\left\|\nabla \Delta \psi\right\|^{2}_{L^{2}}\leq \frac{\mathcal{E}_{0}}{t^{2}}
\]
and we derive the following inequality: for $t>t_{0}$ 
\eqn \label{after t0 2}
\frac{d}{dt} \Big(\frac{\left\|\psi\right\|^{2}_{L^{2}}}{t^{2}}+ \frac{\left\|\nabla \psi\right\|^{2}_{L^{2}}}{t}+\left\|\Delta \psi\right\|^{2}_{L^{2}}\Big) +\mathcal{E}_{0}t^{2}\Big(\frac{\left\|\psi\right\|^{2}_{L^{2}}}{t^{2}}+ \frac{\left\|\nabla \psi\right\|^{2}_{L^{2}}}{t}+\left\|\Delta \psi\right\|^{2}_{L^{2}}\Big)^{2}\leq 0.
\een
Let $\mathcal{I}_{0}=\frac{\left\|\psi(t_{0})\right\|^{2}_{L^{2}}}{t^{2}_{0}}+ \frac{\left\|\nabla \psi(t_{0})\right\|^{2}_{L^{2}}}{t_{0}}+\left\|\Delta \psi(t_{0})\right\|^{2}_{L^{2}}$. By solving (\ref{after t0 2}), we have 
\eqn \label{after t0 3}
\frac{\left\|\psi(t)\right\|^{2}_{L^{2}}}{t^{2}}+ \frac{\left\|\nabla \psi(t)\right\|^{2}_{L^{2}}}{t}+\left\|\Delta \psi(t)\right\|^{2}_{L^{2}}\leq \frac{3\mathcal{I}_{0}}{3+3(t^{3}-t^{3}_{0})\mathcal{I}_{0}}.
\een
Since $ \left\|\psi(t)\right\|_{H^{2}}\leq  \left\|\psi_{0}\right\|_{H^{2}}$ for all $t>0$ by Theorem \ref{GWP}, we obtain 
\[
\left\|\Delta \psi(t)\right\|_{L^{2}}\leq \frac{\mathcal{E}_{0}}{(1+t)^{3/2}}
\]
which complete the proof of Theorem \ref{decay of psi}.

\section{Proof of Theorem \ref{Asymptotics}}
The purpose of this section is to establish the asymptotic behavior of $(\psi, Z)$ as $t\rightarrow \infty$. Let 
\[
\Gamma(t,x)=\frac{1}{4\pi t}e^{-\frac{|x|^{2}}{4t}}
\]
be the two dimensional heat kernel. We first notice that we have the $L^{p}$ estimates of $\Gamma$ in two dimensions: for $1\leq p\leq r\leq \infty$
\begin{equation} \label{decay rate of Gamma}
\begin{split}
&\left\|\Gamma(t) \ast f\right\|_{L^{r}}\leq C(p,r) t^{-\left(\frac{1}{p}-\frac{1}{r}\right)}\|f\|_{L^{p}},\\
&\left\|\nabla \Gamma(t) \ast f\right\|_{L^{r}}\leq C(p,r) t^{-\left(\frac{1}{p}-\frac{1}{r}\right)-\frac{1}{2}}\|f\|_{L^{p}}\\
& \left\|\nabla^{2} \Gamma(t) \ast f\right\|_{L^{r}}\leq C(p,r) t^{-\left(\frac{1}{p}-\frac{1}{r}\right)-1}\|f\|_{L^{p}},
\end{split}
\end{equation}
where $\ast$ is the convolution in the space variables. We also observe that constant multiples of $\Gamma$ are solutions of (\ref{Hall equation in 2D}) because $\Gamma$ and $\Gamma_{t}$ are radial functions and so 
\eqn \label{orthogonality}
\nabla^{\perp}\Gamma\cdot \nabla \Gamma=0, \quad \nabla^{\perp}\Gamma\cdot \nabla \Delta\Gamma=\nabla^{\perp}\Gamma\cdot \nabla \Gamma_{t}=0.
\een

We are now in position to prove Theorem \ref{Asymptotics}. Let $\widetilde{\psi}=\psi-\gamma \Gamma$ and $\widetilde{Z}=Z-\eta\Gamma$, where $\gamma$ and $\eta$ are defined in (\ref{L1 average c}). By using (\ref{orthogonality}), we have 
\eqn \label{tilde equations}
\widetilde{\psi}_{t}-\Delta \widetilde{\psi}=-\nabla^{\perp}Z\cdot \nabla \psi, \quad \widetilde{Z}_{t}-\Delta \widetilde{Z}=-\nabla^{\perp}\psi\cdot \nabla \Delta\psi. 
\een 
So, there are two types of the integral forms of $(\psi,Z)$ from (\ref{Hall equation in 2D}) and (\ref{tilde equations}): 
\begin{subequations} \label{no tilde integral equations}
\begin{align}
& \psi(t)=\Gamma(t) \ast \psi_{0}-\int^{t}_{0} \Gamma(t-s)\ast (\nabla^{\perp}Z\cdot \nabla \psi)(s)ds, \label{no tilde integral equations a}\\
& Z(t)=\Gamma(t) \ast Z_{0}-\int^{t}_{0} \Gamma(t-s)\ast (\nabla^{\perp}\psi\cdot \nabla \Delta\psi)(s)ds \label{no tilde integral equations b}
\end{align}
\end{subequations}
and $\psi=\widetilde{\psi}+\gamma \Gamma$ and $Z=\widetilde{Z}+\eta\Gamma$ with 
\begin{subequations} \label{tilde integral equations}
\begin{align}
& \widetilde{\psi}(t)=\Gamma(t) \ast (\psi_{0}-\gamma\delta_{0})-\int^{t}_{0} \Gamma(t-s)\ast (\nabla^{\perp}Z\cdot \nabla \psi)(s)ds, \label{tilde integral equations a}\\
& \widetilde{Z}(t)=\Gamma(t) \ast (Z_{0}-\eta\delta_{0})-\int^{t}_{0} \Gamma(t-s)\ast (\nabla^{\perp}\psi\cdot \nabla \Delta\psi)(s)ds, \label{tilde integral equations b}
\end{align}
\end{subequations}
where $\delta_{0}$ is the Dirac delta function supported at the origin. Since the time integrals of  (\ref{no tilde integral equations}) and (\ref{tilde integral equations}) are same, the only differences in the asymptotic behaviors are given by the linear parts. In particular, we need (\ref{L1 average c}) to handle the linear part of (\ref{tilde integral equations}). We here estimate $(\widetilde{\psi}, \widetilde{Z})$ which also give the estimation of $(\psi,Z)$.

We now estimate $\widetilde{\psi}$ in $L^{\infty}$ with $\nabla^{\perp}Z\cdot \nabla \psi=\dv \left(\nabla^{\perp}Z \psi\right)$:
\begin{equation*}
\begin{split}
\left\|\widetilde{\psi}(t)\right\|_{L^{\infty}}& \leq \left\|\Gamma(t) \ast (\psi_{0}-\gamma\delta_{0})\right\|_{L^{\infty}}+\int^{\frac{t}{2}}_{0} \left\|\dv \Gamma(t-s)\ast (\nabla^{\perp}Z \psi)(s)\right\|_{L^{\infty}}ds\\
&+ \int^{t}_{\frac{t}{2}} \left\|\dv \Gamma(t-s)\ast (\nabla^{\perp}Z \psi)(s)\right\|_{L^{\infty}}ds=\text{(I)}+\text{(II)}+\text{(III)},
\end{split}
\end{equation*}
We begin with $\text{(I)}$:
\begin{equation*}
\begin{split}
\text{(I)}&=\left\|\int_{\mathbb{R}^{2}}\left(\Gamma(t, x-y)-\Gamma(t,x)\right) \psi_{0}(y)dy \right\|_{L^{\infty}}\\
&=\left\|\int_{\mathbb{R}^{2}}\int^{1}_{0}\nabla \Gamma(t, x-\theta y)\cdot y \psi_{0}(y)d\theta dy \right\|_{L^{\infty}} \leq \left\|\nabla \Gamma(t)\right\|_{L^{\infty}} \left\|\la y\ra \psi\right\|_{L^{1}}\leq \frac{\mathcal{E}_{0}}{t^{\frac{3}{2}}}.
\end{split}
\end{equation*}
To bound $\text{(II)}$, we use Theorem \ref{GWP},  (\ref{L2 bound psi}), and (\ref{decay rate of Gamma}):
\begin{equation*}
\begin{split}
\text{(II)}&\leq C\int^{\frac{t}{2}}_{0}(t-s)^{-\frac{3}{2}} \left\|\nabla^{\perp}Z(s)\psi(s)\right\|_{L^{1}}ds \leq C \int^{\frac{t}{2}}_{0}(t-s)^{-\frac{3}{2}}\left\|\nabla Z(s)\right\|_{L^{2}}\left\|\psi(s)\right\|_{L^{2}}ds \\
&\leq \frac{\mathcal{E}_{0}}{t^{3/2}} \int^{\frac{t}{2}}_{0} \frac{1}{(s+1)\sqrt{s+1}}ds \leq \frac{\mathcal{E}_{0}}{t^{\frac{3}{2}}}.
\end{split}
\end{equation*}
We also bound $\text{(III)}$ by using Theorem \ref{GWP}, Theorem \ref{decay of psi}, (\ref{L2 bound psi}), and (\ref{decay rate of Gamma}):
\begin{equation*}
\begin{split}
\text{(III)}&\leq C\int^{t}_{\frac{t}{2}}(t-s)^{-\frac{5}{6}} \left\|\nabla^{\perp}Z(s)\psi(s)\right\|_{L^{3}}ds \leq C \int^{t}_{\frac{t}{2}}(t-s)^{-\frac{5}{6}}\left\|\nabla Z(s)\right\|_{L^{6}}\left\|\psi(s)\right\|_{L^{6}}ds \\
& \leq C\int^{t}_{\frac{t}{2}}(t-s)^{-\frac{5}{6}}\left\|\nabla Z(s)\right\|^{\frac{1}{3}}_{L^{2}} \left\|\Delta Z(s)\right\|^{\frac{2}{3}}_{L^{2}} \left\|\psi(s)\right\|^{\frac{1}{3}}_{L^{2}} \left\|\nabla \psi(s)\right\|^{\frac{2}{3}}_{L^{2}}ds  \\
&\leq \mathcal{E}_{0} \int^{t}_{\frac{t}{2}}(t-s)^{-\frac{5}{6}} s^{-\frac{1}{6}}s^{-\frac{2}{3}}s^{-\frac{1}{3}}s^{-1}ds \leq \frac{\mathcal{E}_{0}}{t^{2}}.
\end{split}
\end{equation*}
Taking all these bounds into account, we find two types of asymptotic behaviors of $\psi$
\begin{equation*}
\begin{split}
& \psi(t,x)=\gamma\Gamma(t,x)+\widetilde{\psi}(t,x)=\gamma\Gamma(t,x) +O(t^{-3/2}),\\
&\psi(t,x)=\gamma\Gamma(t,x)+\widetilde{\psi}(t,x)=\Gamma(t)\ast \psi_{0} +O(t^{-3/2}).
\end{split}
\end{equation*}

We now derive the same kind of estimates for $\widetilde{Z}$ in $L^{\infty}$ with $\nabla^{\perp}Z\cdot \nabla \Delta \psi=\dv \left(\nabla^{\perp}\psi \Delta \psi\right)$:
\begin{equation*}
\begin{split}
\left\|\widetilde{Z}(t)\right\|_{L^{\infty}}& \leq \left\|\Gamma(t) \ast (Z_{0}-\gamma\delta_{0})\right\|_{L^{\infty}}+\int^{\frac{t}{2}}_{0} \left\|\dv \Gamma(t-s)\ast (\nabla^{\perp}\psi \Delta \psi)(s)\right\|_{L^{\infty}}ds\\
&+ \int^{t}_{\frac{t}{2}} \left\|\dv \Gamma(t-s)\ast (\nabla^{\perp}\psi \Delta \psi)(s)\right\|_{L^{\infty}}ds=\text{(IV)}+\text{(V)}+\text{(VI)}.
\end{split}
\end{equation*}
$\text{(IV)}$ is bounded exactly as $\text{(I)}$:
\[
\text{(IV)} \leq \left\|\nabla \Gamma(t)\right\|_{L^{\infty}} \left\|\la y\ra Z\right\|_{L^{1}}\leq \frac{\mathcal{E}_{0}}{t^{\frac{3}{2}}}.
\]
Before bounding $\text{(V)}$, we rewrite $\dv \Gamma(t-s)\ast (\nabla^{\perp}\psi \Delta \psi)$ as 
\begin{equation*}
\begin{split}
&\dv \Gamma(t-s)\ast (\nabla^{\perp}\psi \Delta \psi)(x)\\
&=-\int \partial_{1}\Gamma(y-x) \Delta \psi(y)\partial_{2}\psi(y)dy+ \int \partial_{2}\Gamma(y-x) \Delta \psi(y)\partial_{1}\psi(y)dy\\
& =\int \partial_{1}\partial_{k}\Gamma(y-x) \partial_{k}\psi(y)\partial_{2}\psi(y)dy-\int \partial_{2}\partial_{k}\Gamma(y-x) \partial_{k}\psi(y)\partial_{1}\psi(y)dy\\
&+ \int \partial_{1}\Gamma(y-x) \partial_{k} \psi(y)\partial_{2}\partial_{k}\psi(y)dy - \int \partial_{2}\Gamma(y-x) \partial_{k} \psi(y)\partial_{1}\partial_{k}\psi(y)dy.
\end{split}
\end{equation*}
We then integrate the last two terms by parts
\begin{equation*}
\begin{split}
&\int \partial_{1}\Gamma(y-x) \partial_{k} \psi(y)\partial_{2}\partial_{k}\psi(y)dy - \int \partial_{2}\Gamma(y-x) \partial_{k} \psi(y)\partial_{1}\partial_{k}\psi(y)dy\\
&=-\int \partial_{1}\partial_{2}\Gamma(y-x) \partial_{k} \psi(y)\partial_{k}\psi(y)dy-\int \partial_{1}\Gamma(y-x) \partial_{2}\partial_{k} \psi(y)\partial_{k}\psi(y)dy\\
&+\int \partial_{1}\partial_{2}\Gamma(y-x) \partial_{k} \psi(y)\partial_{k}\psi(y)dy+\int \partial_{2}\Gamma(y-x) \partial_{1}\partial_{k} \psi(y)\partial_{k}\psi(y)dy
\end{split}
\end{equation*} 
which gives
\[
\int \partial_{1}\Gamma(y-x) \partial_{k} \psi(y)\partial_{2}\partial_{k}\psi(y)dy - \int \partial_{2}\Gamma(y-x) \partial_{k} \psi(y)\partial_{1}\partial_{k}\psi(y)dy=0
\]
and so we obtain 
\[
\begin{split}
\dv \Gamma(t-s)\ast (\nabla^{\perp}\psi \Delta \psi)(x) &=\int \partial_{1}\partial_{k}\Gamma(y-x) \partial_{k}\psi(y)\partial_{2}\psi(y)dy\\
&-\int \partial_{2}\partial_{k}\Gamma(y-x) \partial_{k}\psi(y)\partial_{1}\psi(y)dy.
\end{split}
\]
We now bound $\text{(V)}$ using  Theorem \ref{decay of psi} and  (\ref{decay rate of Gamma}):
\[
\begin{split}
\text{(V)}&\leq C\int^{\frac{t}{2}}_{0}(t-s)^{-2} \left\|\nabla \psi(s)\nabla \psi(s)\right\|_{L^{1}}ds \leq C \int^{\frac{t}{2}}_{0}(t-s)^{-2}\left\|\nabla \psi(s)\right\|^{2}_{L^{2}}ds \\
&\leq \frac{\mathcal{E}_{0}}{t^{2}} \int^{\frac{t}{2}}_{0} \frac{1}{(s+1)^{2}}ds \leq \frac{\mathcal{E}_{0}}{t^{2}}.
\end{split}
\]
We finally bound $\text{(VI)}$ using Theorem \ref{GWP}, Theorem \ref{decay of psi}  and (\ref{decay rate of Gamma}):
\begin{equation*}
\begin{split}
\text{(VI)}&\leq C\int^{t}_{\frac{t}{2}}(t-s)^{-\frac{5}{6}} \left\|\nabla^{\perp}\psi(s)\Delta\psi(s)\right\|_{L^{3}}ds \leq C \int^{t}_{\frac{t}{2}}(t-s)^{-\frac{5}{6}}\left\|\nabla \psi(s)\right\|_{L^{6}}\left\|\Delta \psi(s)\right\|_{L^{6}}ds \\
& \leq C \int^{t}_{\frac{t}{2}}(t-s)^{-\frac{5}{6}}\left\|\nabla \psi(s)\right\|^{\frac{1}{3}}_{L^{2}}\left\|\Delta \psi(s)\right\|_{L^{2}}\left\|\nabla \Delta \psi(s)\right\|^{\frac{2}{3}}_{L^{2}}ds\\
&\leq \mathcal{E}_{0}\int^{t}_{\frac{t}{2}}(t-s)^{-\frac{5}{6}} s^{-\frac{1}{3}}s^{-\frac{3}{2}} s^{-1}ds \leq \frac{\mathcal{E}_{0}}{t^{3}}.
\end{split}
\end{equation*}
Taking all these bounds into account, we also obtain two types of asymptotic behaviors of $Z$: 
\begin{equation*}
\begin{split} 
&Z(t,x)=\eta\Gamma(t,x)+\widetilde{Z}(t,x)=\eta\Gamma(t,x) +O(t^{-\frac{3}{2}}),\\
& Z(t,x)=\eta\Gamma(t,x)+\widetilde{Z}(t,x)=\Gamma(t)\ast Z_{0} +O(t^{-2}).
\end{split}
\end{equation*}

\section{Proof of Theorem \ref{perturbation theorem 1} and Theorem \ref{perturbation theorem 2}} \label{sec:4}
This section is devoted to proving the global existence and the uniqueness of solutions of (\ref{Hall equation in 2D}) around harmonic functions. The analysis here is very close to the one in Section \ref{sec:3.3}, but we will take a different kind of smallness condition, and the existence of harmonic functions requires a bit more computation.

\subsection{Proof of Theorem \ref{perturbation theorem 1}} \label{sec:4.1}
We recall the equations of $\rho$ and $Z$:
\begin{subequations}\label{new perturbation equation 11}
\begin{align}
&\rho_{t}-\Delta \rho=[\rho,Z]+ [\overline{\psi},Z],  \label{new perturbation equation 11 a}\\
&Z_{t} -\Delta Z=[\Delta \rho,\rho]+ [\Delta \rho, \overline{\psi}]  \label{new perturbation equation 11 b}.
\end{align}
\end{subequations}

\subsubsection{\bf A priori estimates} \label{sec:5.1.1}

By (\ref{commutator e}), we have
\eqn \label{perp 1}
\frac{d}{dt}\left(\left\|\nabla \rho\right\|^{2}_{L^{2}}+\left\|Z\right\|^{2}_{L^{2}}\right)+2\left\|\Delta \rho\right\|^{2}_{L^{2}}+2\left\|\nabla Z\right\|^{2}_{L^{2}}=0.
\een

We next multiply (\ref{new perturbation equation 11 a}) by $\Delta^{2} \rho$,  (\ref{new perturbation equation 11 b}) by $-\Delta Z$ and integrate over $\mathbb{R}^{2}$ to get
\begin{equation*}
\begin{split}
&\frac{1}{2}\frac{d}{dt}\left(\left\|\Delta \rho\right\|^{2}_{L^{2}}+\left\|\nabla Z\right\|^{2}_{L^{2}}\right) +\left\|\nabla \Delta \rho\right\|^{2}_{L^{2}}+\left\|\Delta Z\right\|^{2}_{L^{2}} \\
&= \int \Delta^{2} \rho [\rho,Z] -\int\Delta Z[\Delta \rho,\rho]+ \int \Delta^{2} \rho [\overline{\psi},Z] -\int\Delta Z[\Delta \rho,\overline{\psi}]=\text{(I)+(II)+(III)+(IV)}.
\end{split}
\end{equation*}
Treating $\text{(I)+(II)}$ as (\ref{H2 bound 1}) with $1/4$ not $1/2$, we have
\eqn \label{perp new dd}
\text{(I)+(II)}\leq C\left\|\Delta \rho\right\|^{2}_{L^{2}}\left\|\nabla\Delta \rho\right\|^{2}_{L^{2}} +\frac{1}{4}\left\|\Delta Z\right\|^{2}_{L^{2}}.
\een
Since
\begin{equation*}
\begin{split}
\text{(III)}+ \text{(IV)}&=2\int\Delta \rho\left([\overline{\psi}_{x},Z_{x}]+[\overline{\psi}_{y},Z_{y}]\right)\leq C\int\left|\nabla^{2}\overline{\psi}\right| \left|\nabla^{2}Z\right| \left|\nabla^{2}\rho\right|  \\
&\leq C\left\|\nabla^{2}\overline{\psi}\right\|^{2}_{L^{\infty}}\left\|\Delta \rho\right\|^{2}_{L^{2}}+\frac{1}{4}\left\|\Delta Z\right\|^{2}_{L^{2}},
\end{split}
\end{equation*}
we obtain 
\eqn \label{perp 2}
\frac{d}{dt}\left(\left\|\Delta \rho\right\|^{2}_{L^{2}}+\left\|\nabla Z\right\|^{2}_{L^{2}}\right) +\left\|\nabla \Delta \rho\right\|^{2}_{L^{2}}+\left\|\Delta Z\right\|^{2}_{L^{2}} \leq C\left(\left\|\nabla\Delta \rho\right\|^{2}_{L^{2}}+\left\|\nabla^{2}\overline{\psi}\right\|^{2}_{L^{\infty}}\right)\left\|\Delta \rho\right\|^{2}_{L^{2}}.
\een

We finally multiply (\ref{new perturbation equation 11 a}) by $-\Delta^{3} \rho$,  (\ref{new perturbation equation 11 b}) by $\Delta^{2} Z$ and integrate over $\mathbb{R}^{2}$:
\begin{equation*}
\begin{split}
&\frac{1}{2}\frac{d}{dt} \left(\left\|\nabla \Delta \rho\right\|^{2}_{L^{2}}+\left\|\Delta Z\right\|^{2}_{L^{2}}\right) + \left\|\Delta^{2} \rho\right\|^{2}_{L^{2}}+\left\|\nabla \Delta Z\right\|^{2}_{L^{2}} \\
&= -\int \Delta^{3} \rho [\rho,Z] +\int\Delta^{2} Z[\Delta \rho,\rho] -\int \Delta^{3} \rho [\overline{\psi},Z] +\int\Delta^{2} Z[\Delta \rho,\overline{\psi}]=\text{(V)+(VI)+(VII)+(VIII)}.
\end{split}
\end{equation*}
By following (\ref{third term}) and (\ref{fourth term})  in  the proof of Theorem \ref{GWP} with $\frac{1}{2}$ replaced by $\frac{1}{4}$, we have
\eqn \label{perp new dddd}
\text{(V)+(VI)}\leq C\left\|\nabla Z\right\|^{2}_{L^{2}} \left\|\nabla \Delta Z\right\|^{2}_{L^{2}} +C\left\|\Delta \rho\right\|^{2}_{L^{2}} \left\|\Delta^{2} \rho\right\|^{2}_{L^{2}}+\frac{1}{4}\left(\left\|\Delta^{2} \rho\right\|^{2}_{L^{2}}+\left\|\nabla \Delta Z\right\|^{2}_{L^{2}}\right).
\een
The last two terms add up in the following way 
\begin{equation*}
\begin{split}
\text{(VII)}+\text{(VIII)}&=-2\int\Delta \rho\left([\overline{\psi}_{x},\Delta Z_{x}]+[\overline{\psi}_{y},\Delta Z_{y}]\right) -2\int\Delta^{2} \rho\left([\overline{\psi}_{x},Z_{x}]+[\overline{\psi}_{y},Z_{y}]\right)\\
& \leq C\int\left|\nabla^{2}\overline{\psi}\right| \left|\nabla^{3}Z\right| \left|\nabla^{3}\rho\right|+C\int\left|\nabla^{2}\overline{\psi}\right| \left|\nabla^{2}Z\right| \left|\nabla^{4}\rho\right| \\
&\leq C\left\|\nabla^{2}\overline{\psi}\right\|^{2}_{L^{\infty}}\left(\left\|\nabla \Delta \rho\right\|^{2}_{L^{2}}+\left\|\Delta Z\right\|^{2}_{L^{2}}\right)+\frac{1}{4}\left(\left\|\Delta^{2} \rho\right\|^{2}_{L^{2}}+\left\|\nabla \Delta Z\right\|^{2}_{L^{2}}\right).
\end{split}
\end{equation*}
So, we arrive at  
\begin{equation} \label{perp 3}
\begin{split}
&\frac{d}{dt} \left(\left\|\nabla \Delta \rho\right\|^{2}_{L^{2}}+\left\|\Delta Z\right\|^{2}_{L^{2}}\right) + \left\|\Delta^{2} \rho\right\|^{2}_{L^{2}}+\left\|\nabla \Delta Z\right\|^{2}_{L^{2}} \\
&\leq C\left\|\nabla Z\right\|^{2}_{L^{2}} \left\|\nabla \Delta Z\right\|^{2}_{L^{2}} +C\left\|\Delta \rho\right\|^{2}_{L^{2}} \left\|\Delta^{2} \rho\right\|^{2}_{L^{2}}+C\left\|\nabla^{2}\overline{\psi}\right\|^{2}_{L^{\infty}}\left(\left\|\nabla \Delta \rho\right\|^{2}_{L^{2}}+\left\|\Delta Z\right\|^{2}_{L^{2}}\right).
\end{split}
\end{equation}

Let $C_{1}=k\left\|\nabla^{2}\overline{\psi}\right\|^{2}_{L^{\infty}}$ with $k$ large enough which is determined below.  By multiplying (\ref{perp 1}) by $C_{1}^{2}$ and  (\ref{perp 2}) by $C_{1}$ and adding the resulting equations to (\ref{perp 3}), we have
\begin{equation} \label{sum of equations 1}
\begin{split}
&\frac{d}{dt}\left(C_{1}^{2}F_{1}+C_{1}F_{2}+ F_{3}\right)+  2C^{2}_{1}\left(\left\|\Delta \rho\right\|^{2}_{L^{2}}+\left\|\nabla Z\right\|^{2}_{L^{2}}\right)+ C_{1}\left(\left\|\nabla \Delta \rho\right\|^{2}_{L^{2}}+\left\|\Delta Z\right\|^{2}_{L^{2}}\right)+F_{4} \\
&\leq \widehat{C}C_{1}\left\|\nabla^{2}\overline{\psi}\right\|^{2}_{L^{\infty}} \left\|\Delta \rho\right\|^{2}_{L^{2}}+\widehat{C}\left\|\nabla^{2}\overline{\psi}\right\|^{2}_{L^{\infty}}\left(\left\|\nabla \Delta \rho\right\|^{2}_{L^{2}}+\left\|\Delta Z\right\|^{2}_{L^{2}}\right)\\
& +CC_{1}\left\|\Delta \rho\right\|^{2}_{L^{2}}\left\|\nabla \Delta \rho\right\|^{2}_{L^{2}}+C\left\|\nabla Z\right\|^{2}_{L^{2}} \left\|\nabla \Delta Z\right\|^{2}_{L^{2}} +C\left\|\Delta \rho\right\|^{2}_{L^{2}} \left\|\Delta^{2} \rho\right\|^{2}_{L^{2}},
\end{split}
\end{equation}
where we fix two constants by $\widehat{C}$ to determine $k$ and $F_{1}, F_{2}, F_{3}, F_{4}$ are defined in (\ref{several norm 1}). We now choose $k$ such that 
\eqn \label{size k 1}
k>2\widehat{C},  \quad C_{1}>1.
\een
Then, one can easily check that  (\ref{sum of equations 1}) can be reduced to 
\[
\frac{d}{dt}\left(C_{1}^{2}F_{1}+C_{1}F_{2}+ F_{3}\right)+  C_{1}F_{3}+F_{4} \leq C\left(C_{1}^{2}F_{1}+C_{1}F_{2}+ F_{3}\right) \left(C_{1}F_{3}+F_{4}\right).
\]
If $C\epsilon_{2}=C\left(C_{1}^{2}F_{1}(0)+C_{1}F_{2}(0)+ F_{3}(0)\right)<1$, we obtain  the following for all $t>0$:
\[
C_{1}^{2}F_{1}(t)+C_{1}F_{2}(t)+ F_{3}(t)+ (1-C\epsilon_{2})\int^{t}_{0}\left(C_{1}F_{3}(s)+F_{4}(s)\right)ds \leq C_{1}^{2}F_{1}(0)+C_{1}F_{2}(0)+ F_{3}(0).
\]

\subsubsection{ \bf Uniqueness}
Suppose there are two solutions $(\rho_{1}, Z_{1})$ and $(\rho_{2}, Z_{2})$. Let $\rho=\rho_{1}-\rho_{2}$ and $Z=Z_{1}-Z_{2}$. Then, $(\rho,Z)$ satisfies the following equations
\[
\begin{split}
&\rho_{t}-\Delta \rho=[\rho_{1}, Z]+[\rho,Z_{2}]+[\overline{\psi},Z], \\
&Z_{t}-\Delta Z=[\Delta \rho, \rho_{1}]+[\Delta \rho_{2}, \rho]+[\Delta \rho, \overline{\psi}].
\end{split}
\]
Since
\[
-\int \Delta \rho [\overline{\psi},Z] +\int Z [\Delta \rho, \overline{\psi}]=0,
\]
the proof of the uniqueness is  identical to the one in Section \ref{sec:3.1.2}.

\subsection{Proof of Theorem \ref{perturbation theorem 2}}\label{sec:4.2}
We recall the equations of $\psi$ and $\omega$:
\begin{subequations}\label{new perturbation equation 21}
\begin{align}
&\psi_{t}-\Delta \psi=[\psi,\omega]+ [\psi, \overline{Z}],  \label{new perturbation equation 21 a}\\
&\omega_{t} -\Delta \omega=[\Delta \psi,\psi]  \label{new perturbation equation 21 b}.
\end{align}
\end{subequations}

\subsubsection{\bf A priori estimates}
Compared to Theorem \ref{perturbation theorem 1}, we also need the $L^{2}$ bound of $\psi$  to complete the proof of Theorem \ref{perturbation theorem 2}. So, we first have  
\[ \label{new perp 1}
\frac{d}{dt} \left\|\psi\right\|^{2}_{L^{2}} +2\left\|\nabla \psi\right\|^{2}_{L^{2}}=0.
\]

We next multiply (\ref{new perturbation equation 21 a}) by $-\Delta \psi$,  (\ref{new perturbation equation 21 b}) by $\omega$, and integrate over $\mathbb{R}^{2}$. Then,  
\begin{equation*}
\begin{split}
&\frac{1}{2}\frac{d}{dt}\left(\left\|\nabla \psi\right\|^{2}_{L^{2}}+\left\|\omega\right\|^{2}_{L^{2}}\right)+\left\|\Delta \psi\right\|^{2}_{L^{2}}+\left\|\nabla \omega\right\|^{2}_{L^{2}}\\
&=-\int \Delta \psi[\psi,\overline{Z}] \leq C\left\|\nabla \overline{Z}\right\|^{2}_{L^{\infty}}\left\|\nabla \psi\right\|^{2}_{L^{2}}+\frac{1}{2}\left\|\Delta \psi\right\|^{2}_{L^{2}}
\end{split}
\end{equation*}
and so we obtain 
\eqn \label{new perp 2}
\frac{d}{dt}\left(\left\|\nabla \psi\right\|^{2}_{L^{2}}+\left\|\omega\right\|^{2}_{L^{2}}\right)+\left\|\Delta \psi\right\|^{2}_{L^{2}}+\left\|\nabla \omega\right\|^{2}_{L^{2}}\leq C\left\|\nabla \overline{Z}\right\|^{2}_{L^{\infty}}\left\|\nabla \psi\right\|^{2}_{L^{2}}. 
\een

We now multiply (\ref{new perturbation equation 21 a}) by $\Delta^{2} \psi$,  (\ref{new perturbation equation 21 b}) by $-\Delta \omega$, and integrate over $\mathbb{R}^{2}$. Then,  
\begin{equation*}
\begin{split}
&\frac{1}{2}\frac{d}{dt}\left(\left\|\Delta \psi\right\|^{2}_{L^{2}}+\left\|\nabla \omega\right\|^{2}_{L^{2}}\right)+\left\|\nabla \Delta \psi\right\|^{2}_{L^{2}}+\left\|\Delta \omega\right\|^{2}_{L^{2}}\\
&=\int \Delta^{2} \psi[\psi,\omega] -\int \Delta  \omega[\Delta \psi,\psi]+\int \Delta^{2} \psi[\psi,\overline{Z}] =\text{(I)+(II)+(III)}.
\end{split}
\end{equation*}
As (\ref{perp new dd}), we bound $\text{(I)+(II)}$ by  
\[
\text{(I)+(II)}\leq C\left\|\Delta \psi\right\|^{2}_{L^{2}}\left\|\nabla\Delta \psi\right\|^{2}_{L^{2}} +\frac{1}{2}\left\|\Delta \omega\right\|^{2}_{L^{2}}.
\]
And we bound $\text{(III)}$ as 
\[
\text{(III)}=\int \Delta \psi \Delta [\psi,\overline{Z}]=2\int \Delta \psi \left([\psi_{x},\overline{Z}_{x}]+ [\psi_{y},\overline{Z}_{y}]\right) \leq C\left\|\nabla \overline{Z}\right\|^{2}_{L^{\infty}}\left\|\Delta \psi\right\|^{2}_{L^{2}}+\frac{1}{2} \left\|\nabla \Delta \psi\right\|^{2}_{L^{2}}.
\]
So, we obtain 
\eqn \label{new perp 3}
\frac{d}{dt}\left(\left\|\Delta \psi\right\|^{2}_{L^{2}}+\left\|\nabla \omega\right\|^{2}_{L^{2}}\right)+\left\|\nabla \Delta \psi\right\|^{2}_{L^{2}}+\left\|\Delta \omega\right\|^{2}_{L^{2}}\leq C\left(\left\|\nabla^{3}\psi\right\|^{2}_{L^{2}} +\left\|\nabla \overline{Z}\right\|^{2}_{L^{\infty}}\right)\left\|\Delta \psi\right\|^{2}_{L^{2}}. 
\een

 We finally multiply (\ref{new perturbation equation 21 a}) by $-\Delta^{3} \psi$,  (\ref{new perturbation equation 21 b}) by $\Delta^{2} \omega$, and integrate over $\mathbb{R}^{2}$. Then,  
\begin{equation*}
\begin{split}
&\frac{1}{2}\frac{d}{dt}\left(\left\|\nabla \Delta \psi\right\|^{2}_{L^{2}}+\left\|\Delta \omega\right\|^{2}_{L^{2}}\right)+\left\|\Delta^{2} \psi\right\|^{2}_{L^{2}}+\left\|\nabla \Delta \omega\right\|^{2}_{L^{2}}\\
&=-\int \Delta^{3} \psi[\psi,\omega] +\int \Delta^{2}  \omega[\Delta \psi,\psi]-\int \Delta^{3} \psi[\psi,\overline{Z}] =\text{(IV)+(V)+(VI)}.
\end{split}
\end{equation*}
Similar to  (\ref{perp new dddd}), we bound $\text{(IV)+(V)}$ by 
\[
\text{(IV)+(V)}\leq C\left\|\nabla \omega\right\|^{2}_{L^{2}} \left\|\nabla \Delta \omega\right\|^{2}_{L^{2}} +C\left\|\Delta\psi\right\|^{2}_{L^{2}} \left\|\Delta^{2}\psi\right\|^{2}_{L^{2}}+\frac{1}{2}\left\|\Delta^{2}\psi\right\|^{2}_{L^{2}}+ \frac{1}{2}\left\|\nabla\Delta \omega\right\|^{2}_{L^{2}}.
\]
To estimate $\text{(VI)}$, we use
\[
\text{(VI)}=- \int \Delta^{2} \psi  [\Delta\psi,\overline{Z}]-2\int \Delta^{2} \psi [\psi_{x},\overline{Z}_{x}]- 2\int \Delta^{2} \psi [\psi_{y},\overline{Z}_{y}]=\text{(VI)}_{(1)}+\text{(VI)}_{(2)}+\text{(VI)}_{(3)}.
\]
The first term is bounded as above:
\[
\text{(VI)}_{(1)}\leq C\left\|\nabla \overline{Z}\right\|^{2}_{L^{\infty}}\left\|\nabla \Delta \psi\right\|^{2}_{L^{2}}+\frac{1}{4} \left\|\Delta^{2} \psi\right\|^{2}_{L^{2}}.
\]
We next estimate $\text{(VI)}_{(2)}$ as 
\begin{equation*}
\begin{split}
\text{(VI)}_{(2)}&=-2\int \Delta \psi \Delta [\psi_{x},\overline{Z}_{x}]=-2\int \Delta \psi [\Delta  \psi_{x},\overline{Z}_{x}]-4\int \Delta \psi [\psi_{xx},\overline{Z}_{xx}]-4\int \Delta \psi [\psi_{xy},\overline{Z}_{xy}]\\
&=\text{(VI)}_{(2a)}+\text{(VI)}_{(2b)}+\text{(VI)}_{(2c)},
\end{split}
\end{equation*}
where we use the fact that $\overline{Z}_{x}$ is harmonic. $\text{(VI)}_{(2a)}$ is bounded as above:
\[
\text{VI}_{(2a)}\leq C\left\|\nabla \overline{Z}\right\|^{2}_{L^{\infty}}\left\|\nabla \Delta \psi\right\|^{2}_{L^{2}}+\frac{1}{16} \left\|\Delta^{2} \psi\right\|^{2}_{L^{2}}.
\]
By the integration by parts,
\[
\text{(VI)}_{(2b)}=-4\int \overline{Z}_{xx}[\Delta \psi, \psi_{xx}]=4\int \overline{Z}_{x}[\Delta \psi, \psi_{xx}]_{x} \leq C\left\|\nabla \overline{Z}\right\|^{2}_{L^{\infty}}\left\|\nabla \Delta \psi\right\|^{2}_{L^{2}}+\frac{1}{16} \left\|\Delta^{2} \psi\right\|^{2}_{L^{2}}.
\]
Similarly, we obtain
\[
\text{(VI)}_{(2c)}\leq C\left\|\nabla \overline{Z}\right\|^{2}_{L^{\infty}}\left\|\nabla \Delta \psi\right\|^{2}_{L^{2}}+\frac{1}{16} \left\|\Delta^{2} \psi\right\|^{2}_{L^{2}}.
\]
Since $\text{(VI)}_{(2)}$ and $\text{(VI)}_{(2)}$ are of the same form, we obtain
\begin{equation*} \label{new perp 4}
\begin{split}
&\frac{d}{dt}\left(\left\|\nabla \Delta \psi\right\|^{2}_{L^{2}}+\left\|\Delta \omega\right\|^{2}_{L^{2}}\right)+\left\|\Delta^{2} \psi\right\|^{2}_{L^{2}}+\left\|\nabla \Delta \omega\right\|^{2}_{L^{2}} \\
&\leq C\left\|\nabla \omega\right\|^{2}_{L^{2}} \left\|\nabla \Delta \omega\right\|^{2}_{L^{2}} +C\left\|\Delta\psi\right\|^{2}_{L^{2}} \left\|\Delta^{2}\psi\right\|^{2}_{L^{2}}+C\left\|\nabla \overline{Z}\right\|^{2}_{L^{\infty}}\left\|\nabla \Delta \psi\right\|^{2}_{L^{2}}.
\end{split}
\end{equation*}

Let $C_{2}=k\left\|\nabla \overline{Z}\right\|^{2}_{L^{\infty}}$ with $k$ large enough which is determined below.  By following the argument in Section \ref{sec:5.1.1}, we obtain
\begin{equation}\label{sum of equations 2}
\begin{split}
&\frac{d}{dt}\left(C_{2}^{3}\left\|\psi\right\|^{2}_{L^{2}}+ C_{2}^{2}K_{1}+C_{2}K_{2}+ K_{3}\right)+2C^{3}_{2}\left\|\nabla \psi\right\|^{2}_{L^{2}}+C^{2}_{2} \left(\left\|\Delta \psi\right\|^{2}_{L^{2}}+\left\|\nabla \omega\right\|^{2}_{L^{2}}\right)\\
&+C_{2}\left(\left\|\nabla \Delta \psi\right\|^{2}_{L^{2}}+\left\|\Delta \omega\right\|^{2}_{L^{2}}\right)+K_{4} \\
&\leq \widehat{C}C^{2}_{2}\left\|\nabla \overline{Z}\right\|^{2}_{L^{\infty}} \left\|\nabla \psi\right\|^{2}_{L^{2}}+ \widehat{C}C_{2}\left\|\nabla \overline{Z}\right\|^{2}_{L^{\infty}} \left\|\Delta \psi\right\|^{2}_{L^{2}} +\widehat{C}\left\|\nabla \overline{Z}\right\|^{2}_{L^{\infty}} \left\|\nabla \Delta \psi\right\|^{2}_{L^{2}} \\
&+CC_{2}\left\|\Delta \psi\right\|^{2}_{L^{2}}\left\|\nabla \Delta \psi\right\|^{2}_{L^{2}}+C\left\|\Delta \psi\right\|^{2}_{L^{2}}\left\|\Delta^{2} \psi\right\|^{2}_{L^{2}}+C\left\|\nabla \omega \right\|^{2}_{L^{2}}\left\|\nabla \Delta \omega\right\|^{2}_{L^{2}},
\end{split}
\end{equation}
where we fix two constants by $\widehat{C}$ to determine $k$ and $K_{1}, K_{2}, K_{3}, K_{4}$ are defined in (\ref{several norm 2}). We now choose $k$ such that 
\eqn \label{size k 2}
k>2\widehat{C}, \quad C_{2}>\max\{2\widehat{C}, 1\}.
\een
Then, (\ref{sum of equations 2}) can be reduced to 
\begin{equation*}
\begin{split}
&\frac{d}{dt}\left(C_{2}^{3}\left\|\psi\right\|^{2}_{L^{2}}+ C_{2}^{2}K_{1}+C_{2}K_{2}+ K_{3}\right)+ C_{2}^{2}K_{2}+ C_{2}K_{3}+K_{4} \\
&\leq C C_{2}K_{2}K_{3}+CK_{2}K_{4}\leq  C\left(C_{2}^{3}\left\|\psi\right\|^{2}_{L^{2}}+ C_{2}^{2}K_{1}+C_{2}K_{2}+ K_{3}\right) \left(C_{2}^{2}K_{2}+ C_{2}K_{3}+K_{4}\right).
\end{split}
\end{equation*}
If $C\epsilon_{3}=C\left(C_{2}^{3}\left\|\psi(0)\right\|^{2}_{L^{2}}+ C_{2}^{2}K_{1}(0)+C_{2}K_{2}(0)+K_{3}(0)\right)<1$, we obtain 
\begin{equation*} \label{final bound per 2}
\begin{split}
&C_{2}^{3}\left\|\psi(t)\right\|^{2}_{L^{2}}+ C_{2}^{2}K_{1}(t)+C_{2}K_{2}(t)+K_{3}(t)+(1-C\epsilon_{3})\int^{t}_{0}\left(C_{2}^{2}K_{2}(s)+ C_{2}K_{3}(s)+K_{4}(s)\right)ds\\
&\leq C_{2}^{3}\left\|\psi(0)\right\|^{2}_{L^{2}}+ C_{2}^{2}K_{1}(0)+C_{2}K_{2}(0)+ K_{3}(0)
\end{split}
\end{equation*}
for all $t>0$.

\subsubsection{\bf Uniqueness}
Suppose there are two solutions $(\psi_{1}, \omega_{1})$ and $(\psi_{2}, \omega_{2})$. Let $\psi=\psi_{1}-\psi_{2}$ and $\omega=\omega_{1}-\omega_{2}$. Then, $(\psi,\omega)$ satisfies the following equations
\begin{subequations}\label{difference perturbation 1}
\begin{align}
& \psi_{t}-\Delta \psi=[\psi_{1}, \omega]+[\psi,\omega_{2}]+[\psi,\overline{Z}], \label{difference perturbation 1 a}\\
& \omega_{t}-\Delta \omega=[\Delta \psi, \psi_{1}]+[\Delta \psi_{2}, \psi]. \label{difference perturbation 1 b}
\end{align}
\end{subequations}
Compared to (\ref{difference}), there is one extra term $[\psi,\overline{Z}]$. When we multiply (\ref{difference perturbation 1 a}) by $-\Delta \psi$,  (\ref{difference perturbation 1 b}) by $g$, and integrate over $\mathbb{R}^{2}$, this term can be bounded by
\[
-\int \Delta \psi[\psi,\overline{Z}] \leq \left\|\nabla\overline{Z}\right\|_{L^{\infty}}\left\|\nabla \psi\right\|_{L^{2}}\left\|\Delta \psi\right\|_{L^{2}} \leq C\left\|\nabla\overline{Z}\right\|^{2}_{L^{\infty}}\left\|\nabla \psi\right\|^{2}_{L^{2}}+\frac{1}{2}\left\|\Delta \psi\right\|^{2}_{L^{2}}.
\]
By changing the constant from 1 to $\frac{1}{2}$ in front of $\left\|\Delta \psi\right\|^{2}_{L^{2}}$, we can follow Section \ref{sec:3.1.2} for the remaining part to complete the proof of the uniqueness.

\section{Proof of Theorem \ref{LWP Hall MHD} and Theorem \ref{GWP Hall MHD}}\label{sec:7}
In this section, we deal with the $2 \frac{1}{2}$ dimensional Hall MHD. We first recall (\ref{coupled two half}):
\begin{subequations} \label{coupled two half ddddd}
\begin{align}
& \psi_{t}-\Delta \psi =[\psi,Z]-[\psi,\phi], \label{coupled two half ddddd a}\\
& Z_{t}-\Delta Z=[\Delta \psi,\psi]-[Z,\phi]+[W,\psi],\label{coupled two half ddddd b}\\
& W_{t}-\Delta W=-[W,\phi]-[\psi,Z],\label{coupled two half ddddd c}\\
& \Delta \phi_{t}-\Delta^{2} \phi=-[\Delta \phi,\phi]+[\Delta \psi,\psi].\label{coupled two half ddddd d}
\end{align}
\end{subequations}
Proceeding as Theorem \ref{LWP}, we define the following norms: 
\begin{equation} \label{coupled norms}
\begin{split}
 & P(t)=P_{1}(t)+P_{2}(t)+P_{3}(t), \quad Q(t)=P_{2}(t)+P_{3}(t)+P_{4}(t),\\
&P_{1}=\left\|\nabla \psi\right\|^{2}_{L^{2}}+\left\|Z\right\|^{2}_{L^{2}}+\left\|\nabla \phi\right\|^{2}_{L^{2}}+\left\|W\right\|^{2}_{L^{2}}, \\
& P_{2}=\left\|\Delta \psi\right\|^{2}_{L^{2}}+\left\|\nabla Z\right\|^{2}_{L^{2}}+\left\|\Delta \phi\right\|^{2}_{L^{2}}+\left\|\nabla W\right\|^{2}_{L^{2}}, \\
 &P_{3}=\left\|\nabla \Delta \psi\right\|^{2}_{L^{2}}+\left\|\Delta Z\right\|^{2}_{L^{2}}+\left\|\nabla \Delta \phi\right\|^{2}_{L^{2}}+\left\|\Delta W\right\|^{2}_{L^{2}}, \\
 & P_{4}=\left\|\Delta^{2} \psi\right\|^{2}_{L^{2}}+\left\|\nabla \Delta Z\right\|^{2}_{L^{2}}+\left\|\Delta^{2} \phi\right\|^{2}_{L^{2}}+\left\|\nabla \Delta W\right\|^{2}_{L^{2}}.
\end{split}
\end{equation}
We remember  that the Hall term  is the most difficult term to deal with the Hall MHD. So while (\ref{coupled two half ddddd}) and the corresponding spaces (\ref{coupled norms}) look very complicated, the terms resulting from the Hall effect has already been handled in Section \ref{sec:3} and other terms can be treated similarly. Thus,  rather than presenting the proof in great detail, we will present a proof with a few calculations omitted.

\subsection{Proof of Theorem \ref{LWP Hall MHD}}

\subsubsection{\bf A priori estimates}
By multiplying $-\Delta \psi$, $Z$, $W$, $-\phi$ to (\ref{coupled two half ddddd a}), (\ref{coupled two half ddddd b}), (\ref{coupled two half ddddd c}), (\ref{coupled two half ddddd d}), respectively, we first obtain 
\eqn\label{L2 bound Hall MHD}
\frac{1}{2}\frac{d}{dt}P_{1}+P_{2}=0.
\een 

\vspace{1ex}

 By multiplying $\Delta^{2} \psi$, $-\Delta Z$, $-\Delta W$, $\Delta \phi$ to (\ref{coupled two half ddddd a}), (\ref{coupled two half ddddd b}), (\ref{coupled two half ddddd c}), (\ref{coupled two half ddddd d}), respectively, we have  
\begin{equation} \label{P2 P3}
\begin{split}
\frac{1}{2}\frac{d}{dt} P_{2}+P_{3}&=\int \Delta^{2} \psi [\psi,Z] -\int\Delta Z[\Delta \psi,\psi]- \int \Delta^{2}\psi[\psi,\phi] +\int\Delta Z[Z,\phi]\\
&-\int \Delta Z[W,\psi] +\int\Delta W [W,\phi]+\int\Delta W[\psi,Z]+\int \Delta \phi[\Delta \psi,\psi]\\
&=\text{I(a)+I(b)+I(c)+I(d)+I(e)+I(f)+I(g)+I(h)}.
\end{split}
\end{equation}
Treating as (\ref{H2 bound 1}) with $\frac{1}{2}$ replaced with $\frac{1}{8}$, we have 
\[
\text{I(a)+I(b)}\leq  C \left\|\Delta Z\right\|_{L^{2}} \left\|\Delta \psi\right\|_{L^{2}}\left\|\nabla\Delta \psi\right\|_{L^{2}} \leq C\left\|\Delta \psi\right\|^{2}_{L^{2}}\left\|\nabla\Delta \psi\right\|^{2}_{L^{2}} +\frac{1}{8}\left\|\Delta Z\right\|^{2}_{L^{2}}.
\]
After some reduction, $\text{I(c)+I(h)}$ is estimated as 
\begin{equation*} \label{Ich}
\begin{split}
\text{I(c)+I(h)}= &-2\int \Delta \psi \left([\psi_{x},\phi_{x}]+[\psi_{y},\phi_{y}] \right)\leq C\left\|\Delta \phi\right\|_{L^{2}} \left\|\nabla^{2}\psi\right\|^{2}_{L^{4}}  \\
&\leq  C \left\|\Delta \phi\right\|_{L^{2}} \left\|\Delta\psi\right\|_{L^{2}}\left\|\nabla\Delta \psi\right\|_{L^{2}} \leq C\left\|\Delta \phi\right\|^{2}_{L^{2}}\left\|\Delta \psi\right\|^{2}_{L^{2}} +\frac{1}{8}\left\|\nabla \Delta \phi\right\|^{2}_{L^{2}}.
\end{split}
\end{equation*}
By the definition of the commutator and by using after integrating by parts, we have 
\begin{equation*} \label{Id}
\begin{split}
\text{I(d)}+\text{I(f)}&=\int \partial_{k} Z\partial_{k}\nabla^{\perp}\phi \cdot \nabla Z +\int \partial_{k} W\partial_{k}\nabla^{\perp}\phi \cdot \nabla W \\
& \leq C\left\|\nabla Z\right\|^{2}_{L^{2}}\left\|\Delta \phi\right\|^{2}_{L^{2}} +\frac{1}{8}\left\|\Delta Z\right\|^{2}_{L^{2}}+C\left\|\nabla W\right\|^{2}_{L^{2}}\left\|\Delta \phi\right\|^{2}_{L^{2}} +\frac{1}{8}\left\|\Delta W\right\|^{2}_{L^{2}}.
\end{split}
\end{equation*}
We finally bound $\text{I(e)+I(g)}$ as 
\begin{equation*}\label{Ieg}
\begin{split}
\text{I(e)+I(g)}&=\int Z[\Delta \psi,W]+2\int Z[\psi_{x}, W_{x}]+2\int Z[\psi_{y},W_{y}] \leq C \left\|\nabla Z\right\|_{L^{4}} \left\|\nabla W\right\|_{L^{4}}\left\|\Delta \psi\right\|_{L^{2}}\\
& \leq C \left\|\nabla Z\right\|^{2}_{L^{2}} \left\|\Delta \psi\right\|^{2}_{L^{2}}+C \left\|\nabla W\right\|^{2}_{L^{2}} \left\|\Delta \psi\right\|^{2}_{L^{2}}+\frac{1}{8} \left\|\Delta Z\right\|^{2}_{L^{2}} +\frac{1}{8} \left\|\Delta W\right\|^{2}_{L^{2}}.
\end{split}
\end{equation*}
So, we obtain  
\begin{equation}\label{H1 bound}
\begin{split}
\frac{d}{dt} P_{2}+P_{3}&\leq C \left\|\Delta \phi\right\|^{2}_{L^{2}}  \left\|\Delta \psi\right\|^{2}_{L^{2}}+C \left\|\nabla Z\right\|^{2}_{L^{2}} \left\|\Delta \psi\right\|^{2}_{L^{2}}+C \left\|\nabla W\right\|^{2}_{L^{2}} \left\|\Delta \psi\right\|^{2}_{L^{2}}\\
&+C \left\|\Delta \phi\right\|^{2}_{L^{2}} \left\|\nabla W\right\|^{2}_{L^{2}}+C \left\|\Delta \phi\right\|^{2}_{L^{2}}  \left\|\nabla Z\right\|^{2}_{L^{2}} +C\left\|\Delta \psi\right\|^{2}_{L^{2}}\left\|\nabla\Delta \psi\right\|^{2}_{L^{2}}.
\end{split}
\end{equation}

\vspace{1ex}

 By multiplying $-\Delta^{3} \psi$, $\Delta^{2} Z$, $\Delta^{2} W$, $-\Delta^{2} \phi$ to (\ref{coupled two half ddddd a}), (\ref{coupled two half ddddd b}), (\ref{coupled two half ddddd c}), (\ref{coupled two half ddddd d}), respectively, we have  
\begin{equation}\label{P3 P4}
\begin{split}
 \frac{1}{2}\frac{d}{dt}P_{3}+P_{4} &=-\int \Delta^{3} \psi [\psi,Z] +\int\Delta^{2} Z[\Delta \psi,\psi]+ \int \Delta^{3}\psi[\psi,\phi] -\int\Delta^{2} Z[Z,\phi]\\
 &+\int \Delta^{2} Z[W,\psi] -\int\Delta^{2} W [W,\phi]-\int\Delta^{2} W[\psi,Z]-\int \Delta^{2} \phi[\Delta \psi,\psi]\\
&=\text{II(a)+II(b)+II(c)+II(d)+II(e)+II(f)+II(g)+II(h)}.
\end{split}
\end{equation}
By following the computations used for (\ref{H3 bound 2}), we have  
\[
\text{II(a)+II(b)}  \leq C\left\|\Delta Z\right\|^{2}_{L^{2}} \left\|\nabla\Delta \psi\right\|^{2}_{L^{2}} +C\left\|\Delta \psi\right\|^{2}_{L^{2}} \left\|\nabla\Delta \psi\right\|^{4}_{L^{2}}  +\frac{1}{6}\left\|\nabla\Delta Z\right\|^{2}_{L^{2}}+\frac{1}{4} \left\|\Delta^{2}\psi\right\|_{L^{2}}. 
\] 
And we bound $\text{II(c)+II(h)}$ as follows
\begin{equation*} \label{IIch}
\begin{split}
\text{II(c)+II(h)}&=\int \Delta^{2}\psi[\Delta \psi, \phi]+2\int \Delta^{2}\psi[\psi_{x}, \phi_{x}] +2\int \Delta^{2}\psi[\psi_{y}, \phi_{y}]\\
&+ 2\int \Delta\psi[\psi_{x}, \Delta\phi_{x}] +2\int \Delta\psi[\psi_{y}, \Delta\phi_{y}]\\
&=\int \partial_{k}\nabla^{\perp}\phi\cdot \nabla \Delta \psi \partial_{k}\Delta \psi+2\int \Delta^{2}\psi[\psi_{x}, \phi_{x}] +2\int \Delta^{2}\psi[\psi_{y}, \phi_{y}]\\
&+ 2\int \Delta\psi[\psi_{x}, \Delta\phi_{x}] +2\int \Delta\psi[\psi_{y}, \Delta\phi_{y}]\\
&\leq C \left\|\Delta \phi\right\|_{L^{2}}\left\|\nabla \Delta \psi\right\|^{2}_{L^{4}}+C \left\|\Delta \phi\right\|_{L^{4}}\left\| \Delta \phi\right\|_{L^{4}}\left\|\Delta^{2} \psi\right\|_{L^{2}}+C\left\|\Delta \phi\right\|^{2}_{L^{4}} \left\|\Delta^{2} \phi\right\|^{2}_{L^{2}}\\
& \leq C\left\|\Delta \phi\right\|^{2}_{L^{2}} \left\|\nabla \Delta \psi\right\|^{2}_{L^{2}} +C\left\|\Delta \phi\right\|^{2}_{L^{2}}\left\|\nabla\Delta \phi\right\|^{2}_{L^{2}}+C\left\|\Delta \psi\right\|^{2}_{L^{2}} \left\|\nabla \Delta \psi\right\|^{2}_{L^{2}}\\
&+\frac{1}{4} \left\|\Delta^{2} \psi\right\|_{L^{2}} +\frac{1}{2} \left\|\Delta^{2} \phi\right\|_{L^{2}}.
\end{split}
\end{equation*}
We also bound $\text{II(d)}$ as  
\begin{equation*} \label{IId}
\begin{split}
\text{II(d)}&=\int \Delta Z[Z,\Delta \phi]+2\int \Delta Z[Z_{x},\phi_{x}]+2\int \Delta Z[Z_{y},\phi_{y}] \\
& \leq C \left\|\Delta \phi\right\|_{L^{4}}\left\|\nabla Z\right\|_{L^{4}}  \left\|\nabla \Delta Z\right\|_{L^{2}}+C \left\|\Delta \phi\right\|_{L^{2}}\left\|\Delta Z\right\|^{2}_{L^{4}} \\
& \leq C \left\|\Delta \phi\right\|^{2}_{L^{2}}\left\|\nabla\Delta \phi\right\|^{2}_{L^{2}}+C\left\|\nabla Z\right\|^{2}_{L^{2}} \left\|\Delta  Z\right\|^{2}_{L^{2}} +C\left\|\Delta \phi\right\|^{2}_{L^{2}} \left\|\Delta Z\right\|^{2}_{L^{2}} +\frac{1}{6}\left\|\nabla \Delta Z\right\|^{2}_{L^{2}}.
\end{split}
\end{equation*}
Similarly, we bound $\text{II(f)}$ as 
\[
\text{II(f)} \leq C \left\|\Delta \phi\right\|^{2}_{L^{2}}\left\|\nabla\Delta \phi\right\|^{2}_{L^{2}}+C\left\|\nabla W\right\|^{2}_{L^{2}} \left\|\Delta  W\right\|^{2}_{L^{2}} +C\left\|\Delta \phi\right\|^{2}_{L^{2}} \left\|\Delta W\right\|^{2}_{L^{2}} +\frac{1}{4}\left\|\nabla \Delta W\right\|^{2}_{L^{2}}.
\] 
We finally bound $\text{II(e)+II(g)}$ as 
\begin{equation*}\label{IIeg}
\begin{split}
\text{II(e)+II(g)}&=\int Z[\Delta W,\Delta \psi]+2\int Z[\Delta W_{x}, \psi_{x}]+2\int Z[\Delta W_{y},\psi_{y}]+\int \Delta Z[W,\Delta \psi]\\
&+2\int \Delta Z[W_{x}, \psi_{x}] +2\int \Delta Z[W_{y}, \psi_{y}]\\
& \leq C \left\|\nabla Z\right\|_{L^{4}}\left\|\Delta \psi\right\|_{L^{4}}\left\|\nabla \Delta W\right\|_{L^{2}} +C \left\|\nabla W\right\|_{L^{4}}\left\|\Delta \psi\right\|_{L^{4}}\left\|\nabla \Delta Z\right\|_{L^{2}}\\
&\leq C \left\|\nabla Z\right\|^{2}_{L^{2}}\left\|\Delta Z\right\|^{2}_{L^{2}}+C \left\|\nabla W\right\|^{2}_{L^{2}}\left\|\Delta W\right\|^{2}_{L^{2}}+ C\left\|\Delta \psi\right\|^{2}_{L^{2}} \left\|\nabla \Delta \psi\right\|^{2}_{L^{2}}  \\
&+\frac{1}{4} \left\|\nabla \Delta W\right\|^{2}_{L^{2}} +\frac{1}{6} \left\|\nabla \Delta Z\right\|_{L^{2}}.
\end{split}
\end{equation*}
Collecting all the bounds, we derive  
\begin{equation}\label{H2 bound}
\begin{split}
\frac{d}{dt} P_{3}+P_{4}&\leq C\left\|\Delta \phi\right\|^{2}_{L^{2}} \left\|\nabla \Delta \psi\right\|^{2}_{L^{2}} +C\left\|\Delta \phi\right\|^{2}_{L^{2}}\left\|\nabla \Delta \phi\right\|^{2}_{L^{2}}+C\left\|\Delta \psi\right\|^{2}_{L^{2}} \left\|\nabla \Delta \psi\right\|^{2}_{L^{2}}\\
&+C\left\|\Delta Z\right\|^{2}_{L^{2}} \left\|\nabla\Delta \psi\right\|^{2}_{L^{2}}  +C \left\|\nabla W\right\|^{2}_{L^{2}}\left\|\Delta W\right\|^{2}_{L^{2}} +C\left\|\Delta \phi\right\|^{2}_{L^{2}} \left\|\Delta W\right\|^{2}_{L^{2}}\\
&+C\left\|\nabla Z\right\|^{2}_{L^{2}} \left\|\Delta  Z\right\|^{2}_{L^{2}} +C\left\|\Delta \phi\right\|^{2}_{L^{2}} \left\|\Delta Z\right\|^{2}_{L^{2}}+C\left\|\Delta \psi\right\|^{2}_{L^{2}} \left\|\nabla\Delta \psi\right\|^{4}_{L^{2}}.
\end{split}
\end{equation}

\vspace{1ex}

By (\ref{L2 bound Hall MHD}), (\ref{H1 bound}), and (\ref{H2 bound}), 
\eqn \label{full H3 d}
\left(1+P(t)\right)'+Q(t) \leq C P^{2}(t)+CP^{3}(t)\leq C \left(1+P(t)\right)^{3}
\een
from which we deduce
\eqn \label{full H3 dd}
P(t)\leq \sqrt{\frac{(1+P(0))^{2}}{1-2Ct(1+P(0))^{2}}}-1 \quad \text{for all} \ t\leq T_{\ast}<\frac{1}{2C (1+P(0))^{2}}.
\een
Integrating (\ref{full H3 d}) and using (\ref{full H3 dd}), we finally derive 
\eqn \label{last bound Hall MHD}
P(t) +\int^{t}_{0}Q(s)ds<\infty, \quad 0<t<T_{\ast}.
\een

\subsubsection{\bf Uniqueness}
Suppose there are two solutions $(\psi_{1}, Z_{1}, \phi_{1}, W_{1})$ and $(\psi_{2}, Z_{2}, \phi_{2}, W_{2})$. Let $\psi=\psi_{1}-\psi_{2}$, $Z=Z_{1}-Z_{2}$,  $\phi=\phi_{1}-\phi_{2}$ and $W=W_{1}-W_{2}$. Then, $(\psi,Z, \phi, W)$ satisfies the following equations
\[
\begin{split}
& \psi_{t}-\Delta \psi =[\psi_{1},Z]+[\psi,Z_{2}]-[\psi_{1},\phi] -[\psi,\phi_{2}], \\
& Z_{t}-\Delta Z=[\Delta \psi,\psi_{1}]+[\Delta \psi_{2},\psi]-[Z_{1},\phi]-[Z,\phi_{2}]+[W_{1},\psi] +[W,\psi_{2}],\\
& W_{t}-\Delta W=-[W_{1},\phi] -[W,\phi_{2}]-[\psi_{1},Z] -[\psi,Z_{2}],\\
& \Delta \phi_{t}-\Delta^{2} \phi=-[\Delta \phi_{1},\phi] -[\Delta \phi,\phi_{2}]+[\Delta \psi_{1},\psi] +[\Delta \psi,\psi_{2}].
\end{split}
\]
From this, we obtain 
\begin{equation*}
\begin{split}
&\frac{1}{2}\frac{d}{dt}\left(\left\|\nabla \psi\right\|^{2}_{L^{2}}+\left\|Z\right\|^{2}_{L^{2}}+\left\|\nabla \phi\right\|^{2}_{L^{2}}+\left\|W\right\|^{2}_{L^{2}}\right)+\left\|\Delta \psi\right\|^{2}_{L^{2}}+\left\|\nabla Z\right\|^{2}_{L^{2}}+\left\|\Delta \phi\right\|^{2}_{L^{2}}+\left\|\nabla W\right\|^{2}_{L^{2}}\\
&=-\int \Delta \psi[\psi_{1},Z]-\int \Delta \psi[\psi,Z_{2}]+\int Z[\Delta \psi,\psi_{1}]+\int Z[\Delta \psi_{2},\psi]\\
&-\int Z[Z_{1},\phi]+\int Z[W_{1},\psi] +\int Z[W,\psi_{2}]-\int W[W_{1},\phi] -\int W[\psi_{1},Z] -\int W[\psi,Z_{2}]\\
&+\int \Delta \psi[\psi_{1},\phi] +\int \Delta \psi[\psi,\phi_{2}]+\int \phi[\Delta \phi,\phi_{2}] -\int \phi[\Delta \psi,\psi_{2}]-\int \phi[\Delta \psi_{1},\psi].
\end{split}
\end{equation*}
The first line on the right-hand side is bounded as (\ref{difference uniqueness})
\[
C\left(\left\|\nabla Z_{2}\right\|^{2}_{L^{2}} \left\|\Delta Z_{2}\right\|^{2}_{L^{2}} +\left\|\Delta \psi_{2}\right\|^{2}_{L^{2}} \left\|\nabla \Delta \psi_{2}\right\|^{2}_{L^{2}} \right)\left\|\nabla \psi\right\|^{2}_{L^{2}} +\frac{1}{3}\left\|\Delta \psi\right\|^{2}_{L^{2}}+ \left\|\nabla Z\right\|^{2}_{L^{2}}
\] 
and the second line is bounded by 
\[
\begin{split}
&C\left\|\nabla Z_{1}\right\|_{L^{\infty}}\left\|Z\right\|_{L^{2}} \left\|\nabla \phi\right\|_{L^{2}} +C\left\|\nabla W_{1}\right\|_{L^{\infty}}\left\|Z\right\|_{L^{2}} \left\|\nabla \psi\right\|_{L^{2}} +C\left\|\nabla \psi_{2}\right\|_{L^{\infty}}\left\|Z\right\|_{L^{2}} \left\|\nabla W\right\|_{L^{2}}\\
& +C\left\|\nabla W_{1}\right\|_{L^{\infty}}\left\|W\right\|_{L^{2}} \left\|\nabla \phi\right\|_{L^{2}} +C\left\|\nabla \psi_{1}\right\|_{L^{\infty}}\left\|W\right\|_{L^{2}} \left\|\nabla Z\right\|_{L^{2}} +C\left\|\nabla Z_{2}\right\|_{L^{\infty}}\left\|W\right\|_{L^{2}} \left\|\nabla \psi\right\|_{L^{2}}.
\end{split}
\]
The third line except for the last one is bounded by 
\[
\begin{split}
&C\left\|\nabla \psi_{1}\right\|_{L^{\infty}}\left\|\nabla \phi\right\|_{L^{2}} \left\|\Delta \psi\right\|_{L^{2}} +C\left\|\nabla \psi_{2}\right\|_{L^{\infty}}\left\|\nabla \psi\right\|_{L^{2}} \left\|\Delta \psi\right\|_{L^{2}} +C\left\|\nabla \phi_{1}\right\|_{L^{\infty}}\left\|\nabla \phi\right\|_{L^{2}} \left\|\Delta \phi\right\|_{L^{2}}  \\
&+C\left\|\nabla \psi_{2}\right\|_{L^{\infty}}\left\|\nabla \phi\right\|_{L^{2}} \left\|\Delta \phi\right\|_{L^{2}} +C\left\|\Delta \psi_{1}\right\|_{L^{2}}\left\|\nabla \phi\right\|_{L^{4}} \left\|\nabla \psi\right\|_{L^{4}}\\
&\leq C\left\|\nabla \psi_{1}\right\|^{2}_{L^{\infty}}\left\|\nabla \phi\right\|^{2}_{L^{2}} +C\left\|\nabla \psi_{2}\right\|^{2}_{L^{\infty}}\left\|\nabla \psi\right\|^{2}_{L^{2}}+C\left\|\nabla \phi_{1}\right\|^{2}_{L^{\infty}}\left\|\nabla \phi\right\|^{2}_{L^{2}}+C\left\|\nabla \psi_{2}\right\|^{2}_{L^{\infty}}\left\|\nabla \phi\right\|^{2}_{L^{2}} \\
&+C \left\|\Delta \psi_{1}\right\|^{2}_{L^{2}}\left\|\nabla \phi\right\|^{2}_{L^{2}} +C \left\|\Delta \psi_{1}\right\|^{2}_{L^{2}}\left\|\nabla \psi\right\|^{2}_{L^{2}}+\frac{2}{3}\left\|\Delta \psi\right\|^{2}_{L^{2}}\left\|\Delta \phi\right\|^{2}_{L^{2}}.
\end{split}
\]
Thus, we arrive at the the following inequality:
\[
\frac{d}{dt}\left(\left\|\nabla \psi\right\|^{2}_{L^{2}}+\left\|Z\right\|^{2}_{L^{2}}+\left\|\nabla \phi\right\|^{2}_{L^{2}}+\left\|W\right\|^{2}_{L^{2}}\right)\leq C\mathcal{I}\left(\left\|\nabla \psi\right\|^{2}_{L^{2}}+\left\|Z\right\|^{2}_{L^{2}}+\left\|\nabla \phi\right\|^{2}_{L^{2}}+\left\|W\right\|^{2}_{L^{2}}\right),
\]
where
\[
\begin{split}
\mathcal{I}&=\left\|\nabla Z_{2}\right\|^{2}_{L^{2}} \left\|\Delta Z_{2}\right\|^{2}_{L^{2}} +\left\|\Delta \psi_{2}\right\|^{2}_{L^{2}} \left\|\nabla \Delta \psi_{2}\right\|^{2}_{L^{2}}+\left\|\nabla Z_{1}\right\|_{L^{\infty}}+\left\|\nabla W_{1}\right\|_{L^{\infty}}+\left\|\nabla \psi_{1}\right\|_{L^{\infty}}\\
&+\left\|\nabla \psi_{2}\right\|_{L^{\infty}}+\left\|\nabla Z_{2}\right\|_{L^{\infty}}+\left\|\nabla \psi_{1}\right\|^{2}_{L^{\infty}} +\left\|\nabla \psi_{2}\right\|^{2}_{L^{\infty}}+\left\|\nabla \phi_{1}\right\|^{2}_{L^{\infty}}+\left\|\Delta \psi_{1}\right\|^{2}_{L^{2}}.
\end{split}
\]
Since $\mathcal{I}$ is integrable in time by (\ref{last bound Hall MHD}), we conclude the uniqueness of solutions.

\subsubsection{\bf Blow-up criterion}
To find a blow-up criterion, we first use (\ref{blowup first}) to bound $\text{I(a)+I(b)}$ as 
\[
\begin{split}
&\left|\int \Delta \psi \left([\psi_{x},Z_{x}]+[\psi_{y},Z_{y}] \right)\right| \leq C\left\|\nabla Z\right\|^{q}_{L^{p}} \left\|\Delta\psi\right\|^{2}_{L^{2}} +\epsilon\left\|\nabla\Delta\psi\right\|^{2}_{L^{2}}, \quad \frac{1}{p}+\frac{1}{q}=\frac{1}{2}. 
\end{split}
\]
So, we can rewrite (\ref{H1 bound}) as 
\begin{equation*}
\begin{split}
\frac{d}{dt}P_{2}  +P_{3}&\leq C \left\|\Delta \phi\right\|^{2}_{L^{2}}  \left\|\Delta \psi\right\|^{2}_{L^{2}}+C \left\|\nabla Z\right\|^{2}_{L^{2}} \left\|\Delta \psi\right\|^{2}_{L^{2}}+C \left\|\nabla W\right\|^{2}_{L^{2}} \left\|\Delta \psi\right\|^{2}_{L^{2}}+C \left\|\Delta \phi\right\|^{2}_{L^{2}} \left\|\nabla W\right\|^{2}_{L^{2}}\\
&+C \left\|\Delta \phi\right\|^{2}_{L^{2}}  \left\|\nabla Z\right\|^{2}_{L^{2}}  +C\left\|\nabla Z\right\|^{q}_{L^{p}} \left\|\Delta\psi\right\|^{2}_{L^{2}}.
\end{split}
\end{equation*}
Integrating this in time with the aid of  (\ref{L2 bound Hall MHD}), we obtain  
\begin{equation}\label{H1 bound New}
\begin{split}
& \left\|\Delta \psi(t)\right\|^{2}_{L^{2}}+\left\|\nabla Z(t)\right\|^{2}_{L^{2}}+\left\|\Delta \phi(t)\right\|^{2}_{L^{2}}+\left\|\nabla W(t)\right\|^{2}_{L^{2}}\\
&+\int^{t}_{0}\left(\left\|\nabla \Delta \psi(s)\right\|^{2}_{L^{2}}+\left\|\Delta Z(s)\right\|^{2}_{L^{2}}+\left\|\nabla \Delta \phi(s)\right\|^{2}_{L^{2}}+\left\|\Delta W(s)\right\|^{2}_{L^{2}}\right) ds\\
&\leq C\mathcal{E}_{0}\exp\left[C\int^{t}_{0}\left\|\nabla Z(s)\right\|^{q}_{L^{p}}ds\right].
\end{split}
\end{equation}
Using the idea in Section \ref{sec:3.1.3} to bound (\ref{H3 bound 2}) by using (\ref{blowup 2}), we bound (\ref{H2 bound}) as follows 
\[
\left\|\nabla\Delta \psi(t)\right\|^{2}_{L^{2}}+\left\|\Delta Z(t)\right\|^{2}_{L^{2}}+\left\|\nabla\Delta \phi(t)\right\|^{2}_{L^{2}}+\left\|\Delta W(t)\right\|^{2}_{L^{2}}\leq \mathcal{E}_{0}\exp \exp\left[CB(t)+CB^{2}(t)\right],
\] 
where $B(t)$ is defined in (\ref{Blowup Hall equation}). This completes the proof of Theorem \ref{LWP Hall MHD}.

\subsection{Proof of Theorem \ref{GWP Hall MHD}}
To prove Theorem \ref{GWP Hall MHD}, we need to bound the quantities used in the proof of Theorem \ref{LWP Hall MHD} in different ways. We first bound each term on the right-hand side of (\ref{P2 P3}) as follow
\begin{equation*}\label{Ibc New}
\begin{split}
\text{I(a)+I(b)}&\leq  C \left\|\Delta Z\right\|_{L^{2}} \left\|\Delta \psi\right\|_{L^{2}}\left\|\nabla\Delta \psi\right\|_{L^{2}} \leq C\left\|\Delta \psi\right\|^{2}_{L^{2}}\left\|\nabla\Delta \psi\right\|^{2}_{L^{2}} +\frac{1}{8}\left\|\Delta Z\right\|^{2}_{L^{2}},\\
\text{I(c)+I(h)}&\leq C \left\|\nabla \phi\right\|^{2}_{L^{2}} \left\|\nabla \Delta \phi\right\|^{2}_{L^{2}}+C\left\|\nabla \psi\right\|^{2}_{L^{2}} \left\|\nabla \Delta \psi\right\|^{2}_{L^{2}}+\frac{1}{8}\left\|\nabla \Delta \psi\right\|^{2}_{L^{2}},\\
\text{I(e)+I(g)}& \leq C \left\|Z\right\|^{2}_{L^{2}}\left\|\Delta Z\right\|^{2}_{L^{2}}+C \left\|W\right\|^{2}_{L^{2}}\left\|\Delta W\right\|^{2}_{L^{2}}+C \left\|\nabla \psi\right\|^{2}_{L^{2}}\left\|\nabla \Delta \psi\right\|^{2}_{L^{2}}\\
&+\frac{1}{8} \left\|\Delta Z\right\|^{2}_{L^{2}} +\frac{1}{8} \left\|\Delta W\right\|^{2}_{L^{2}},\\
\text{I(d)+I(f)}& \leq C   \left\|\nabla \phi\right\|^{2}_{L^{2}}\left\|\nabla \Delta \phi\right\|^{2}_{L^{2}} +C\left\|W\right\|^{2}_{L^{2}}\left\|\Delta W\right\|^{2}_{L^{2}}+C \left\|Z\right\|^{2}_{L^{2}}\left\|\Delta Z\right\|^{2}_{L^{2}}\\
&+\frac{1}{8} \left\|\Delta W\right\|^{2}_{L^{2}} +\frac{1}{8} \left\|\Delta Z\right\|^{2}_{L^{2}}.
\end{split}
\end{equation*}
So, we can rewrite (\ref{H1 bound}) as 
\eqn \label{H1 bound New 2}
\frac{d}{dt} P_{2}+P_{3}\leq C \left(\left\|\nabla \phi\right\|^{2}_{L^{2}} +\left\|\nabla \psi\right\|^{2}_{L^{2}}+\left\|Z\right\|^{2}_{L^{2}} +\left\|W\right\|^{2}_{L^{2}} +\left\|\Delta \psi\right\|^{2}_{L^{2}} \right)P_{3}.
\een 
By (\ref{L2 bound Hall MHD}) and (\ref{H1 bound New 2}),
\eqn \label{P1 plus P2}
\frac{d}{dt}(P_{1}+P_{2})+P_{2}+P_{3}\leq C \left(P_{1}+P_{2}\right)\left(P_{2}+P_{3}\right).
\een
Let $\epsilon_{4}= \left\|\nabla \psi_{0}\right\|^{2}_{H^{1}}+\left\|Z_{0}\right\|^{2}_{H^{1}}+\left\|\nabla \phi_{0}\right\|^{2}_{H^{1}}+\left\|W_{0}\right\|^{2}_{H^{1}}$. If $\epsilon_{4}$ is sufficiently small such as $C\epsilon_{4}<1$, (\ref{P1 plus P2}) implies the following inequality for all $t>0$:
\eqn \label{P1 plus P2 integral}
P_{1}(t)+P_{2}(t)+(1-C\epsilon_{4})\int^{t}_{0}\left(P_{2}(s)+P_{3}(s)\right)ds\leq \epsilon_{4}.
\een

We also bound the right-hand side of (\ref{P3 P4}) as follows
\begin{equation*} \label{IIab New}
\begin{split}
\text{II(a)+II(b)} & \leq C\left\|\nabla Z\right\|^{2}_{L^{2}}\left\|\nabla \Delta Z\right\|^{2}_{L^{2}} +C\left\|\Delta \psi\right\|^{2}_{L^{2}}\left\|\Delta^{2} \psi\right\|^{2}_{L^{2}} +C\left\|\Delta \psi\right\|^{3}_{L^{2}} \left\|\Delta^{2} \psi\right\|^{2}_{L^{2}}  \\
&+\frac{1}{6}\left\|\nabla\Delta Z\right\|^{2}_{L^{2}}+\frac{1}{4} \left\|\Delta^{2}\psi\right\|_{L^{2}},\\
\text{II(c)+II(h)}&  \leq C\left\|\nabla \phi\right\|^{2}_{L^{2}}  \left\|\Delta^{2} \psi\right\|^{2}_{L^{2}} +C  \left\|\nabla \psi\right\|^{2}_{L^{2}} \left\|\Delta^{2}\phi\right\|^{2}_{L^{2}}  +C\left\|\nabla \phi\right\|^{2}_{L^{2}}\left\|\Delta^{2}\phi\right\|^{2}_{L^{2}}\\
& +C\left\|\nabla \psi\right\|^{2}_{L^{2}} \left\|\Delta^{2} \psi\right\|^{2}_{L^{2}}+\frac{1}{4} \left\|\Delta^{2} \psi\right\|_{L^{2}},\\
\text{II(e)+II(g)}&\leq C \left\| Z\right\|^{2}_{L^{2}}\left\|\nabla \Delta Z\right\|^{2}_{L^{2}}+C \left\|W\right\|^{2}_{L^{2}}\left\|\nabla \Delta W\right\|^{2}_{L^{2}}+ C\left\|\nabla \psi\right\|^{2}_{L^{2}} \left\|\Delta^{2} \psi\right\|^{2}_{L^{2}}  \\
&+\frac{1}{4} \left\|\nabla \Delta W\right\|^{2}_{L^{2}} +\frac{1}{6} \left\|\nabla \Delta Z\right\|_{L^{2}},\\
\text{II(d)+II(f)}& \leq C \left\|\nabla \phi\right\|^{2}_{L^{2}}\left\|\Delta^{2}\phi\right\|^{2}_{L^{2}}+C\left\| W\right\|^{2}_{L^{2}} \left\|\nabla \Delta  W\right\|^{2}_{L^{2}} +C\left\|Z\right\|^{2}_{L^{2}} \left\|\nabla \Delta  Z\right\|^{2}_{L^{2}}\\
&+C\left\|\nabla \phi\right\|^{2}_{L^{2}}\left\|\nabla \Delta W\right\|^{2}_{L^{2}}+C \left\|W\right\|^{2}_{L^{2}}\left\|\Delta^{2}\phi\right\|^{2}_{L^{2}}   +C\left\|\nabla \phi\right\|^{2}_{L^{2}}\left\|\nabla \Delta Z\right\|^{2}_{L^{2}}\\
& +C \left\|Z\right\|^{2}_{L^{2}}\left\|\Delta^{2}\phi\right\|^{2}_{L^{2}} +\frac{1}{4}\left\|\nabla \Delta W\right\|^{2}_{L^{2}}+\frac{1}{6}\left\|\nabla \Delta Z\right\|^{2}_{L^{2}}.
\end{split}
\end{equation*}
So, we rewrite (\ref{H2 bound}) as 
\[
\frac{d}{dt} P_{3}+P_{4} \leq C\left(\left\|Z\right\|^{2}_{L^{2}}+\left\|\nabla Z\right\|^{2}_{L^{2}}+\left\|W\right\|^{2}_{L^{2}}+\left\|\nabla \psi \right\|^{2}_{L^{2}}+\left\|\nabla \phi\right\|^{2}_{L^{2}}+\left\|\Delta \psi\right\|^{2}_{L^{2}} \right)P_{4}.
\] 
By (\ref{P1 plus P2 integral}), we derive the following for all $t>0$:
\eqn \label{P3 integral}
P_{3}(t)+\left(1-C\epsilon_{4}\right)\int^{t}_{0}P_{4}(s)ds\leq P_{3}(0).
\een

Combining (\ref{L2 bound Hall MHD}), (\ref{P1 plus P2 integral}), (\ref{P3 integral}), we conclude that 
\[
P(t)+\left(1-C\epsilon_{0}\right)\int^{t}_{0}Q(s)ds\leq P(0)
\]
for all $t>0$. This completes the proof of Theorem \ref{GWP Hall MHD}.

\section*{Acknowledgments}
H.B. was supported by NRF-2018R1D1A1B07049015.  

K.K. was supported by NRF-2019R1A2C1084685 and NRF-20151009350.

\bibliographystyle{amsplain}

\bibliographystyle{amsplain}

\end{document}